\documentclass[12pt]{amsart}
\usepackage{amsmath}
\usepackage{amssymb}
\usepackage{amsfonts}
\usepackage{amsthm}
\usepackage{verbatim}
\usepackage{bookman}
 \DeclareMathOperator{\Spf}{Spf}
 \DeclareMathOperator{\Deff}{Deff}
 \DeclareMathOperator{\ext}{Ext}
 \DeclareMathOperator{\NCI}{NCI}
 
\def\inv{^{-1}}
\DeclareMathOperator{\codim}{codim}
\DeclareMathOperator{\xc}{XC}
\DeclareMathOperator{\del}{\partial}
\DeclareMathOperator{\Ext}{Ext} \DeclareMathOperator{\Hom}{Hom}
\def\refp #1.{(\ref{#1})}

\newcommand{\bsl}{\begin{slide}{}}

\newcommand{\h}{\mathcal{H}}

\newcommand{\st}{\ ^{\tilde{}}}

\newcommand{\A}{\mathcal{A}}

\newcommand{\X}{\mathcal{X}}
\newcommand{\eps}{\epsilon}

\newcommand{\qis}{\sim_{\mathrm{qis}}}
\newcommand{\ba}{{\beta\alpha}}

\newcommand{\Cal}[1]{\mathcal #1}
\newcommand{\act}[2]{\langle #1 , #2\rangle}

\def\sbr #1.{^{[#1]}}
\def\sfl #1.{^{\lfloor #1\rfloor}}
\newcommand{\nred}{\ {^{\mathrm{red}}\mathcal{N}}}
\newcommand{\tred}{{\ ^{\mathrm{red}}\mathcal{T}}}
\newcommand{\ttt}{\tilde{\mathcal{T}}}
\newcommand{\tn}{\tilde{\mathcal{N}}}

\newcommand {\txp}{\t_X(P)}
\newcommand\jsh{\rm J^\sharp}
\newcommand\gb{\mathfrak g_\bullet}
\newcommand\half{\frac{1}{2}}
\newcommand\red{{\mathrm red}}

\def\what{\widehat}
\def\inv{^{-1}}
\def\?{{\bf{??}}}
\def\dgla{differential graded Lie algebra\ }

\newcommand{\F}{\mathcal F}
\newcommand\ul[1]{\underline #1}

\def\H{\Cal H}
\def\HH{\mathbb H}

\def\A{\Bbb A}
\def\X{\Cal X}

\def\F{\Cal F}

\def\C{\mathbb C}
\def\P{\mathbb P}

\def\Z{\mathbb Z}

\def\Der{\text{\rm Der} }
\def\Spec{\text{\rm Spec} }
\def\ad{\text{\rm ad}}

\def\O{\mathcal O}

\def\Sym{\textrm{Sym}}
\def\id{\text{id}}

\def\g{\mathfrak g}

\def\gl{\mathfrak{gl}}

\def\h{\mathfrak h}
\def\rs{{\rho\sigma}}
\def\ab{{\alpha\beta}}
\def\bg{{\beta\gamma}}
\def\ga{{\gamma\alpha}}
\def\abg{{\alpha\beta\gamma}}
\def\m{\mathfrak m}

\def\hom{\mathfrak {hom}}

\def\1/2{\frac{1}{2}}

\def\D{\mathfrak D}
\def\U{\mathfrak U}
\def\I{\mathcal{ I}}

\def\2{{[2]}}
\def\l{\ell}
\def\nl{\newline}
\def\n{\mathcal{N}}
\def\t{\mathcal{T}}
\def\hom{\mathcal{H}\mathit{om}}
\def\RHom{\text{\rm R}\mathcal H\mathit om}
\def\That{\hat{T}}
\def\<{\langle}
\def\>{\rangle}
\def\sela{SELA\ }
\def\endo{\mathcal{E}\mathit{nd}}

\def\beq{\begin{equation}}
\def\eeq{\end{equation}}
\newcommand{\mex}[1]{\begin{example}\begin{rm}#1
\end{rm}
\end{example}}
\newcommand \beql[2] {\begin{equation}\label{#1}#2\end{equation}}
\def\eex{\end{rm}\end{example}}
\newcommand\newsection[1]{\section{#1}\setcounter{equation}{0}}

\newtheorem{thm}{Theorem}[section]
\newtheorem{cor}[thm]{Corollary}

\newtheorem{lem}[thm]{Lemma}

\newtheorem{defn}[thm]{Definition}
\theoremstyle{remark}
\newtheorem{rem}[thm]{Remark}
\newtheorem{example}[thm]{Example}

\begin{document}
\title{ Jacobi-Bernoulli cohomology \\ and deformations of schemes and
maps}

\author{Ziv Ran}
\thanks{
Research supported in part by NSA grant
 H98230-05-1-0063 and by the Mittag-Leffler
Institute 06-07 Moduli Spaces program}
\address{University of California,
 Riverside}
\email{ziv.ran@ucr.edu}
\begin{abstract} We introduce a notion of
Jacobi-Bernoulli
cohomology associated to a semi-simplicial Lie algebra (SELA).
For
an algebraic scheme $X$ over $\C$, we construct a tangent SELA
$\t_X$ and show that the Jacobi-Bernoulli cohomology of $\t_X$
is
related to infinitesimal deformations of $X$.
\end{abstract}

\maketitle \setcounter{section}{-1} \section{Overview} The
'usual'
deformation theory, e.g. of complex structures, in the manner
of
Kodaira-Spencer-Grothendieck (cf. e.g. \cite{K2, Sernesi} and
references
therein), is commonly couched in terms of a differential graded
Lie
algebra or dgla $\g$. It can be viewed, as in \cite{cid}, as
studying the canonically defined \emph{deformation ring} $R(\g)$,
fashioned from
the
\emph{Jacobi cohomology}, i.e. the cohomology of the Jacobi
complex
associated to $\g$. This setting is somewhat restrictive, e.g.
it is
not broad enough to accomodate such naturally occurring
deformation
problems as embedded deformations of a submanifold $X$ in a
fixed
ambient space $Y$. In \cite{atom} we introduced the notion of
\emph{Lie atom} (essentially, Lie pair) and an associated
\emph{Jacobi-Bernoulli complex} as an extension of that of dgla
and
its Jacobi complex, one that is broad enough to handle embedded
deformations and a number of other problems besides.\par A
purpose
of this paper is to establish the familiar notion of (dg)
\emph{semi-simplicial} Lie algebra (\sela) as an appropriately
general and convenient setting for deformation theory. As a
first
approximation, one can think of \sela as a structure like that
of
the \v{C}ech complex of a sheaf of Lie algebras on a
topological
space $X$ with respect to some open covering of $X$. Not only
is
\sela a broad generalization of Lie atom, it is broad enough,
as
we
show, to encompass deformations of (arbitrarily singular)
algebraic
schemes (over $\C$).\par
To compute the deformation theory of a
\sela $\gb$ we introduce a complex that we call the
Jacobi-Bernoulli
complex of $\gb$, though a more proper attribution would be to
Jacobi-Bernoulli-Baker-Campbell-Hausdorff. In a nutshell, the
point
of this complex is that it transforms a gluing condition from
nonabelian coycle condition to ordinary (additive) cocycle
condition
via the multilinearity of the groups making up the complex. A
typical gluing condition looks like
\begin{equation*}
\Psi_{\alpha\beta}\Psi_{\beta\gamma}\Psi_{\gamma\alpha}=1
\end{equation*} with
$\Psi_{\alpha\beta}\in\exp(\g_{\alpha\beta})$,
where $\g_\ab$ may be thought of as the component of our \sela
$\gb$ having to do with gluing over $U_\alpha\cap U_\beta$.
This
condition can be transformed as follows. Write
$$\Psi_{\alpha\beta}=\exp(\psi_{\alpha\beta})$$ etc. Now the
BCH
formula gives a formal expression
$$\exp(X)\exp(Y)\exp(Z)=\sum W_{i,j,k}(X,Y,Z)$$ where
$W_{i,j,k}(X,Y,Z)$ is a homogeneous ad-polynomial of tridegree
$i,j,k$ (`BCH polynomial'), which can be viewed as a linear map
$$w_{ijk}:\Sym^i(\g_\ab)\otimes\Sym^j(\g_\bg)\otimes\Sym^k(\g_\ga)\to
g_\abg$$ Then the above gluing condition on the $\Psi_\bullet$
 becomes the additive condition
that
$$w_{i,j,k}(\psi_\ab^i\otimes\psi_\bg^j\otimes\psi_\ga^k)=0,\
\forall
i,j,k.$$ Now our Jacobi-Bernoulli complex $J(\gb)$ for the
\sela
$\gb$ is essentially designed so as to encompass the various
BCH
polynomials $w_{ijk}$. It is a comultiplicative complex whose
groups
essentially constitute the symmetric algebra on $\gb$ and whose
maps
are essentially derived from the $w_{ijk}$ by the requirement
of
comultiplicativity. The dual of the cohomology of $J(\gb)$
yields
the deformation ring associated to the \sela $\gb$.\par As
mentioned
above, our other main result here is that the deformation
theory
of
an algebraic scheme over $\C$ can be expressed in terms of a
\sela.
Unsurprisingly, this is done via an affine covering. Thus the
first
step is to associate a dgla to a closed embedding $$X\to P$$
where
$P$ is an affine (or for that matter, projective) space. We
call
this the \emph{tangent} dgla to $X$ and denote it $\t_X(P)$. In
a
nutshell, $\t_X(P)$ is defined as the mapping cone of a map
that
we
construct
$$T_P\otimes A_X\to N_{X/P}$$ where $N_{X/P}$ is the normal
atom
to
$X$ in $P$ as in \cite{atom}. That is, $\t_X(P)$ is represented
by
the mapping cone of a map of free modules representing
$T_P\otimes A_X$ and $N_{X/P}$. We will show $\t_X(P)$ admits a
dgla
structure, a dgla action on $A_X$, as well as $A_X$-module
structure. Up to a certain type of 'weak equivalence', the dgla
$\t_X(P)$ depends only on the isomorphism class of $X$ and not
on
the embedding in $P$.
\par The partial independence on the embedding is good enough
to
enable us to associate a global \sela $\t_{X\bullet}$ for an
arbitrary algebraic scheme $X$ defined in terms of, but up to
weak
equivalence independent of,  an affine covering $X_\alpha$ and
embeddings of each $X_\alpha$ in an affine space $P_\alpha$:
e.g.
$$\t_{X,\alpha}=\t_{X_\alpha}(P_\alpha),$$
$$\t_{\X, \ab}=\t_{\X_\alpha\cap X_{\beta}}(P_\alpha\times
P_\beta)$$ etc. Global deformations of $X$ then amount to a
collections of deformations of each $X_\alpha$, given via
Kodaira-Spencer theory by a suitable element
$\phi_\alpha\in\t_{X,\alpha}^1$, plus a collection of gluing
data
$\psi_\ab\in\t_{X,\ab}^0,$ and the necessary compatibilities
are
readily expressed as a cocycle condition in the
Jacobi-Bernoulli
complex $J(\t_{X\bullet}).$
\begin{rem}
At the time of writing (Sep '07), this preprint is not yet in final form; in particular, it has portions developed
at different times which have yet to be fully synchronized.
Comments or questions from readers would be most welcome.
\end{rem}

\section{Semi-Simplicial Lie algebras and Jacobi-Bernoulli
complex }
\subsection{\sela} Our
notion of \sela is essentially the dual of the portion of the
usual
notion of simplicial Lie algebra involving only the face maps
without degeneracy maps. Let $A$ be a totally ordered index-set.
A
\emph{simplex} in $A$ is a finite nonempty subset $S\subset A$,
while a \emph{biplex} is a pair $(S_1\subset S_2)$ of simplices
with
$|S_1|+1=|S_2|$; similarly for \emph{triplex} $(S_1\subset S_2\subset
S_3)$
etc. The \emph{sign} $\epsilon(S_1,S_2)$ of a biplex
$(S_1,S_2)$
is
defined by the condition that
$$\epsilon((0,...,\hat{p},...,n), (0,...,n))=(-1)^{n-p}.$$ By a
\emph{semi simplicial Lie algebra} (SELA) $\gb$ on $A$ we shall mean
the
assignment for each simplex $S$ on $A$ of a Lie algebra $\g_S$,
and
for each biplex $(S_1, S_2)$ of a map ('coface' or
'restriction')
$$r_{S_1,S_2}:\g_{S_1}\to\g_{S_2}$$ such that
$\epsilon(S_1,S_2)r(S_1,S_2)$ is a Lie homomorphism and such
that
for each $S_1\subset S_3$ with $|S_1|+2=|S_3|$, we have
\begin{equation}
\label{complex}\sum_{\stackrel{\rm triplex} {(S_1\subset
S_2\subset
S_3)}}r_{(S_2,S_3)} r_{(S_1,S_2)}=0.\end{equation} The identity
(\ref{complex}) implies that we may assemble the $\g_S$ into a
complex $K^.(\gb)$ where
$$K^i(\gb)=\bigoplus\limits_{|S|=i+1}\g_S$$ and differential
constructed from the various $r_{(S_1,S_2)}$.
\begin{example}{\begin{rm} If $\g$ is a sheaf of Lie algebras
on
a
topological space $X$, and $(U_\alpha)$ is an open covering of
$X$,
there is a \v{C}ech \sela $$\gb: S\mapsto
\g(\bigcap\limits_{\alpha\in
S}U_\alpha).$$ The standard complex $K^.(\gb)$ is
the
\v{C}ech complex $\check{C}(\g,(U_\alpha))$; it  plays a
fundamental role in the study of $\g$-deformations.
\end{rm}}\end{example}
\mex{If $\g\to\h$ is a Lie pair (more generally, a Lie atom,
cf.
\cite{atom}), we get a \sela $\g_.$ on the index-set $(01)$ with $\g_0=\g,
\g_1=0,
\g_{01}=\h.$\vskip 1cm
\begin{picture}(60,20)(20,10)
\put(195, 35){$\h$}
\put(175,30){$\g_\bullet$\line(1,0){30}$_\bullet
0$}\end{picture}\newline
 The deformation-theoretic
significance of $\g_.$ is like that of the Lie atom $(\g,\h)$,
viz.
$\g$-deformations together with an $\h$-trivialization.\par An
obvious generalization would be to take a pair of maps
$\g_1\to\h,
\g_2\to\h$ (e.g. twice the same map), which corresponds to
pairs\
($\g_1$- deformation, $\g_2$-deformation)\ such that
the induced
$\h$-deformations are equivalent as such.}
\subsection{Bernoulli numbers and Baker-Campbell-Hausdorff}
Let $\g$ be a
nilpotent Lie algebra. For an element $X\in\g$ we consider the
formal exponential $\exp(X)$ as an element of the formal
enveloping
algbera $\U(\g).$ Then we can write
\begin{equation}\label{beta}\exp(X)\exp(Y)=\exp(\beta(X,Y))
\end{equation} where $\beta$ is a certain
bracket-polynomial in $X,Y$, known as the
Baker-Campbell-Hausdorff
or BCH polynomial. We denote by $\beta_{i,j}, \beta_i$ the
portion
of $\beta$ in bidegree $i,j$ (resp. total degree $i$). Note
that
each $\beta_{i,j}(X,Y)$ will be a linear combination  of
(noncommutative) $\ad$ monomials with a total of $i$ many $X$'s
and
$j$ many $Y$'s. We write such a monomial in the form
\beq\label{admon1}\ad_S(X^iY^j)=
\ad(T_1)\circ...\circ\ad(T_{i+j-1})(T_{i+j})\eeq
where $S\subset [1,i+j]$ is a subset of cardinality $i$ and
$T_k=X$
(resp. $T_k=Y$) iff $k\in S$ (resp. $k\not\in S$). We denote by
\beq\ad_S(X_1,...,X_i, Y_1,...,Y_j)\eeq the analogous
function, 
obtained by replacing the $x$th occurrence of $X$ (resp. $y$th
occurrence of $Y$) by an $X_x$ (resp. $Y_y$) and by
$\ad_S^\Sym(X_1,...,X_i,Y_1,...,Y_j)$ the corresponding
symmetrized
version, i.e.
\beq\label{admon3}\ad_S^\Sym(X_1,...,X_i,Y_1,...,Y_j)=\sum\limits_{
\stackrel{\pi\in\mathfrak S_i}{\rho\in\mathfrak
S_j}}\frac{1}{i!j!}
\ad_S(X_{\pi(1)},...,X_{\pi(i)},Y_{\rho(1)},...,Y_{\rho(j)}).\eeq\par
We will compute $\beta$, following \cite{vara}, \S2.15 (where
Varadarajan attributes the argument to lectures of Bargmann
that
follow original papers by Baker and Hausdorff). Set
\begin{equation}\label{cd}D(x)=\frac{e^x-1}{x}, C(x)=1/D(x).
\end{equation} Thus, $C(x)$ is the generating function for the
Bernoulli
numbers $B_n$, i.e.
$$C(x)=1+\sum\limits_{n=1}^\infty \frac{B_n}{n!}x^n
=\sum\limits_{n=0}^\infty C_nx^n.$$ Now the reader can easily
check
that for any derivation $\partial$ we have
$$\partial\exp(U)\exp(-U)=D(\ad (U))(\partial U),
\exp(-U)\partial\exp(U)=D(-\ad(U))(\partial U).$$ Now
differentiate
(\ref{beta}) with respect to $X$ and multiply both sides by\nl
$\exp(-\beta(X,Y)).$ This yields (where $\partial_X$ is the
unique
derivation taking $X$ to $X$ and $Y$ to 0)\par
\lefteqn{ X=\partial_X(\exp(X))\exp(-X)=}
$$~~~~~~~~~~~~~~~~~~
{=\partial_X(\exp(\beta(X,Y)))\exp(-\beta(X,Y))
=D(\ad(\beta(X,Y)))(\partial_X\beta(X,Y)).}$$
Thus
\begin{equation}\label{betax}
\partial_X\beta(X,Y)=C(\ad(\beta(X,Y)))(X).
\end{equation} Similarly,
\begin{equation}\label{betay}\partial_Y\beta(X,Y)=
C(-\ad(\beta(X,Y)))(Y).
\end{equation} Starting from $\beta_0=0$, the formulas
(\ref{betax}),(\ref{betay}) clearly determine $\beta$. For
example,
clearly $\beta_{0,*}(X,Y)=Y$, therefore it follows that
\begin{equation}\beta_{1,*}(X,Y)=
C(\ad(Y))(X)=X+\frac{1}{2}[X,Y]+\frac{1}{12}\ad(Y)^2(X)+...
\end{equation} \par We
shall require the obvious extension of this set-up to the
trivariate
case. Thus define a function $\beta(X,Y,Z)$ (NB same letter as
for
the bivariate version) by
\beq\exp(X)\exp(Y)\exp(Z)=\exp(\beta(X,Y,Z))\eeq and let
$\beta_{i,j,k}$ denote its portion in tridegree $(i,j,k)$. Note
that
$$\beta(X,Y,Z)=\beta(\beta(X,Y),Z).$$
\subsection{Jacobi-Bernoulli complex}
Let $\gb$ be a \sela. For simplicity, we shall assume $\gb$ is
2- dimensional, in the sense that $\g_S=0$ for any simplex $S$
of
dimension $>2$; for our applications to deformation theory,
this
is
not a significant restriction. We will also assume that $\gb$
is
\emph{strongly nilpotent} in the sense that it is an algebra
over a
commutative ring $R$ such that $\g_S^{\otimes N}=0$ for all
simplices $S$ and some integer $N$ independent of $S$, with all
tensor products over $R$. This condition obviously depends only
on
the $S$-module structure of $\gb$ and not on its Lie bracket.
We
are
going to define a filtered complex $J=\jsh_m(\gb).$ The groups
$J^j$
can be defined succinctly as
$$J^j=(\Sym^.(K^.(\gb)[1]))^j$$ where $\Sym^.$ is understood in
the
signed or graded sense, alternating on odd terms, and
$K^.(\gb)[1]$
is the standard complex on $\gb$ shifted left once (which is a
complex in degrees $-1,0,1$) . The increasing filtration $F.$
is
by
`number of multiplicands'.  More concretely,\beq
J^{j,k}=\bigoplus\limits_
{\stackrel{-\sum\limits_i\l_i+\sum\limits_in_i=j}
{\sum\limits_i\l_i+\sum\limits_im_i+\sum\limits_in_i=
k}}\bigotimes\limits^i\bigwedge\limits^{\l_i}\g_{\alpha_i}
\otimes
\bigotimes\limits^i\Sym^{m_i}\g_{\alpha_i\beta_i}\otimes
\bigotimes\limits^i\bigwedge\limits^{n_i}\g_{\alpha_i\beta_i\gamma_i}
\eeq \beq F_mJ^j=\bigoplus\limits_{k\leq m}J^{j,k},\eeq \beq
J^j=F_\infty J^j=F_NJ^j.\eeq To define the differential $d$ on
$J^.$, we proceed in steps. Let $\alpha<\beta<\gamma$ be
indices
and
recall that we are identifying $\g_{\gamma\alpha}$ with
$\g_{\alpha\gamma}$.\begin{itemize}\item The differential is
defined
so that the obvious inclusion \beq K^.(\gb)[1])=F_1J^.\to
J^.\eeq is
a map of complexes.
\item The
component\begin{equation*} \Sym^i\g_{\gamma\alpha}\otimes
\Sym^j\g_{\ab}\otimes\Sym^\g_{\beta\gamma}\to
\g_{\alpha\beta\gamma}\end{equation*} is given
by\begin{equation}
 {X^iY^jZ^k\mapsto\beta_{i,j,k}(X,Y,Z)}.\end{equation}
\item The component
\begin{equation*}\nonumber
\g_\alpha\otimes\Sym^i\g_{\alpha\beta}\otimes
\Sym^n\g_{\beta\gamma}
\to\Sym^{i-t+1}\g_{\alpha\beta}\otimes\Sym^n\g_{\beta\gamma},
0\leq
t\leq i\end{equation*} is given by\begin{equation} X\otimes
Y^i\otimes Z^n\mapsto C_tY^{i-t}\ad(Y)^t(X)\otimes
Z^n\end{equation}(where $C_t$ is the normalized Bernoulli
coefficient).\item All componets are defined subject to the
'derivation rule', which means commutativity
of the following diagram
$$\forall i=i_1+i_2, j=j_1+j_2,\ \  i_1,i_2,j_1, j_2>0:$$
\begin{equation}\label{derivation}\begin{matrix}
&J_i&\to&J_{i_i}\otimes J_{i_2}&\\
\partial_{i,j}&\downarrow&&\downarrow&\partial_{i_1,j_1}\otimes
1_{i_2,j_2} + 1_{i_1,j_1}\otimes\partial_{i_2,j_2}
\\ &J_j&\to&J_{j_1}\otimes J_{j_2}&
\end{matrix}\end{equation} in which $J_i=\bigoplus\limits_j
J^{j,i}$
denotes the group
direct summand of $J$ in multiplicative degree $i$ and
$\partial_{i,j}$ denotes the component of $\partial$ going from
$J_i$ to $J_j$, and $1_{i,j}$ is the identity if $i=j$ and zero
otherwise; and of course the horizontal arrows are the
appropriate components of the comultiplication map.
E.g. the component
$$\g_\alpha\otimes\Sym^i\g_{\beta\gamma}\otimes\Sym^k\g_{\beta\gamma}
\to\g_{\gamma\alpha}\otimes
\Sym^i\g_{\ab}\otimes\Sym^k\g_{\beta\gamma}$$ is extended in
the
obvious way from the given differential
$\g_\alpha\to\g_{\gamma\alpha}$ which comes with the \sela
data.
\item Components not defined via the
above rules are set equal to 0. In particular, the component
$$\g_\alpha\otimes\g_{\alpha\beta\gamma}\to
\g_{\alpha\beta\gamma}$$
is zero.
\end{itemize} The following result summarizes the main
properties of
the Jacobi-Bernoulli complex $J$ associated to a \sela (not
least,
that it is a complex!). It is in part, but not entirely, a
direct
extension of the analogous result for Lie atoms given in
\cite{atom}.
\begin{thm}\label{tang sela}
\begin{enumerate}\item $(J^., F.)$ is a functor from the
category of SELAs over $S$ to that of comultiplicative,
cocummutative and coassociative filtered complexes over
$S$.\item
The filtration $F_.$ is compatible with the comultiplication
and
has
associated graded
$$F_i/F_{i-1}=\bigwedge\limits^iK^.(\gb).$$\item
 The quasi-isomorphism class of $J(\gb)$ depends only on the
quasi-isomorphism class of $\gb$ as \sela.\end{enumerate}
\end{thm}\begin{proof} As in the proof of \cite{atom}, Thm
1.2.1,
the main issue is to prove $J$ is a complex, i.e. $d^2=0$. And
again
as in \cite{atom}, it suffices, in light of the derivation
rule,
to
prove the vanishing of the components of $d^2$ that land in
$F_1$,
i.e. that have just one multiplicative factor. Among those, the
proof that these components of $d^2$ vanish on terms of degree
$\leq
-2$, i.e. involving $\bigwedge\limits^i\g_\alpha, i\geq 2$, is
similar to the case of the JB complex considered in
\cite{atom}.
The
essential new case, not considered in \cite{atom}, is the
vanishing
of the $F_1$- components of $d^2$ on terms of degree $-1$, i.e.
terms of the from
$$X\otimes Y^i\otimes Z^n\in
\g_\alpha\otimes\Sym^i\g_{\alpha\beta}\otimes
\Sym^n\g_{\beta\gamma}.$$ For such a term, what needs to be
shown is
the vanishing of the component of $d^2$ of it in
$\g_{\alpha\beta\gamma}$. Thus, we need to prove that\beq d^2(
X\otimes Y^i\otimes Z^n)_{\g_{\alpha\beta\gamma}}=0.\eeq Now
this
component gets contributions via the various components of $d(
X\otimes Y^i\otimes Z^n)$ and those contributions come in  two
kinds:\begin{itemize}\item Via
$g_{\gamma\alpha}\otimes\Sym^i\g_{\alpha\beta}\otimes
\Sym^n\g_{\beta\gamma}$, we get $-\beta_{1,i,n}(X,Y,Z)$. This
comes
from $\beta(\beta(X,Y),Z)$, but is only affected by the terms
in
$\beta(X,Y)$ of degree $\leq 1$ in $X$, i.e. by
$$U=Y+C(\ad(Y))(X)=Y+\sum\limits_{t=0}^\infty C_t\ad(Y)^t(X).$$
This
contribution is obtained by taking $-\beta_{i+1-t,n}(U,Z)$ and
replacing each monomial $$\ad_S(U^{i+1-t}Z^n)$$ (cf.
(\ref{admon1}-\ref{admon3}) ) by
$$(i+1-t)\ad_S^\Sym(C_t(\ad(Y)^t(X))Y^{i-t}Z^n)$$
where
and finally summing as $t$ ranges from 0 to $i$.
\item Via
$\Sym^{i+1-t}\g_{\alpha\beta}\otimes\Sym^n\g_{\beta\gamma}$,
for
each $0\leq t\leq i$, we get a contribution equal to the
expression
obtained by taking $\beta_{i+1-t,n}(W,Z)$, replacing each
monomial
$$\ad_S(W^{i+1-t}Z^n)$$ {\ \ \ \textrm {by}\ \ \ }
$$(i+1-t)\ad_S^\Sym(C_t\ad(Y)^t(X)Y^{i-t}Z^n).$$
\end{itemize} Thus, the sum total of all contributions to
$d^2(X\otimes Y^i\otimes Z^n)_{\g_{\alpha\beta\gamma}}$ is
zero.
\end{proof} The local ring
$$(R(\gb),\m_{R(\gb)})=(\C\oplus\HH^0(J(\gb))^*,\HH^0(J(\gb))^*)$$ is called the
\emph{deformation
ring} of $\gb$. Its (more meaningful) $\m$-adic completion is denoted $\hat{R}(\gb)$ and called the \emph{formal deformation
ring} of $\gb$. The formal spectrum $\Spf(\hat{R}(\gb)$ is called the \emph{formal deformation space} of $\gb$ and denoted $\Deff(\gb)$.
\begin{cor}Assume that $K^.(\gb)$ is acyclic in
nonpositive  degrees. Then there is a second-quadrant spectral sequence
with $E_1$ term
\beq\label{specseq}E_1^{p,q}=\Sym^{-q-2p}\H^1(\gb)
\otimes\Sym^{q+p}\H^2(\gb)\eeq
whose abutment has degree 0 part $\bigoplus E^{p,-p}_\infty$ equalt to
\beq\label{specseq-e_infty} E^{p,-p}_\infty=\bigoplus {\text{gr}}^{-p}_{F_.}\H^0(J(\gb)).\eeq
\end{cor}
\begin{proof} Identical to that of \cite{atom}, Cor. 1.2.2.\end{proof}
\begin{cor}\label{dim-ineq} If $\gb$ is acyclic in nonpositive degrees
and $h^1(\gb)$ and $h^2(\gb)$ are finite, then
we have a lower bound for the Krull dimension of
its formal deformation ring:
$$\dim\hat{R}(\gb)\geq h^1(\gb)-h^2(\gb).$$\end{cor}
\begin{proof} Cf. \cite{atom}, Cor. 1.2.3.\end{proof}
\subsection{Special multiplicative cocycles} Let $(S,\m_S)$ be
a
local artin $\C$-algebra and let $\gb^.=\gb$ be a dg-\sela
(i.e.
for each simplex $T$,
$\g_T$ is a dgla and the coface maps are $\pm$dg
homomorphisms).
A
special class of (multiplicative) cocycles for the
Jacobi-Bernoulli
complex $$J^.(\gb\otimes\m_S)\subset J^.(\gb)\otimes\m_S$$ can
be
constructed as follows. Suppose $$\phi_\bullet\in
K^0(\gb)^1\otimes\m_S=\bigoplus\limits_\rho\g_\rho^1\otimes\m_S,$$
$$\psi_\bullet\in
K^1(\gb)^0\otimes\m_S=\bigoplus_\rs\g_\rs^0\otimes\m_S$$ are
such
that, $\forall \rho,\sigma,\tau$,
\begin{eqnarray}\label{special1}
\partial\phi_\rho&=&-\half [\phi_\rho,\phi_\rho],\\
\label{special2}
\partial\psi_\rs&=&-C(\ad(\psi_\rs))(\phi_\rho)+
C(-\ad(\psi_{\rho\sigma}))(\phi_\sigma),\\
\label{special3}\beta(\psi_\rs, \psi_{\sigma\tau},
\psi_{\tau\rho})&=&0
 .\end{eqnarray}

Then let $$\epsilon(\phi_\bullet, \psi_\bullet)\in
J^0(\gb)\otimes\m_S$$ be the element with components
$$\psi_\rs\in
\g_\rs\otimes\m_S,$$$$\phi_\rho\in\g_\rho\otimes\m_S, $$ and
generally \begin{equation}\label{epsilon}
\bigwedge\limits^r\phi_\rho\otimes\bigotimes
(\psi_{\rho_i\sigma_i})^{n_i}\in
\bigwedge\limits^r(\g^1_\rho\otimes\m_S)
\otimes\bigotimes\Sym^{n_i}(\g_{\rho_i\sigma_i}\otimes\m_S),
 r, n \geq 0.\end{equation} We call
$\epsilon(\phi_\bullet, \phi_\bullet)$ a \emph{special
multiplicative cocycle with coefficients in $S$}.
\begin{lem}\begin{enumerate}\item Assuming
(\ref{special1}-\ref{special3}),
the cochain $\epsilon(\phi_\bullet, \psi_\bullet)$
defined in \refp{epsilon}. is a 0-cocycle for
$J(\gb\otimes\m_S)$
and the associated map
$$^t[\epsilon(\phi_\bullet, \psi_\bullet)]\in\mathrm{Hom}(
R(\gb),S)$$
is a local  ring homomorphism.\item Assume $K^.(\gb)$
is acyclic in negative degrees. Then given an artin local
algebra
$S$, there is a
natural bijection
between the set of cohomology classes
$$[\epsilon(\phi_\bullet, \psi_\bullet)]
\in H^0(J^.(\gb\otimes\m_S))$$
  of special multiplicative cocycles with
coefficients in $S$ and the set of local ring homomorphisms
$R(\gb)\to
S.$\end{enumerate}
\end{lem}\begin{proof}(i) Because of the derivation property
of the differentials in the complex $J$,
it will suffice to prove the vanishing of the portion of the
differential
$\partial \epsilon(\phi_\bullet, \psi_\bullet)$ in $F_1J$
(i.e. having a single tensor factor). That portion
had just three components: in $K^0(\gb)\otimes\m_S, K^1(\gb)^1$
and $K^2(\gb)^0$, and the vanishing of these yields,
respectively eq. (\ref{special1}-\ref{special3}).\par
Now given that $\epsilon $ is a cocycle, the fact that
$^t[\epsilon]$ is multiplicative follows from the fact,
clear by construction, that
$\epsilon= \epsilon(\phi_\bullet, \psi_\bullet)$ is compatible
with the comultiplication on $J^.=J^.(\gb\otimes\m_S)$, in the
sense that the image of $\epsilon$ in $J^.\otimes J^.$ is just
$\epsilon\otimes\epsilon$.
\par (ii)
Note that because $K^.(\gb)$ is acyclic in negative degrees,
the natural map $$H^0(J^.(\gb\otimes\m_S))\to
H^0(J^.(\gb))\otimes\m_S$$ is injective, so any class
$\gamma\in H^0(J^.(\gb\otimes\m_S))$ is determined by
the associated map $^t\gamma:R(\gb)\to S$ . Hence, our
assertion
is a matter of showing that for any local
$\C$-algebra homomorphism $f:R(\gb)\to S$, the corresponding
cohomology class\nl $\eta\in H^0(J^.(\gb))\otimes\m_S$ can be
represented in the form $[\epsilon(\phi_\bullet,
\psi_\bullet)]$, where
$$(\phi_\bullet,\psi_\bullet)\in
(K^0(\gb)^1\otimes K^1(\gb)^0)\otimes\m_S$$
is just the component of (a representative of) $\eta$
of multiplicative degree 1, i.e.
in $F_1(J)$, where $J=J^.(\gb)\otimes\m_S$.
We may assume $f:R_k(\gb)\to S$ where
$k$ is the exponent of $S$ and, by induction, that the result
holds
for $S$ of exponent $<k$. Set
$\epsilon=\epsilon(\phi_\bullet,
\psi_\bullet)$. By induction, $\epsilon$ is a cocycle modulo
the
socle $I={\mathrm{ann}}_S(\m_S)$ of $S$ and $\eta\equiv
[\epsilon]\mod I$. Moreover,
because $\m^2_{R_k(\gb)}$ and $\m^2_S$ are quotients,
respectively, of $\m_{R_{k-1}(\gb)}$ and $\m_{S/I}$, the images
of $\eta$ and $[\epsilon]$ in $(J/F_1(J))\otimes\m_S$ (whic are
automatically in $(J/F_1(J))\otimes\m_S^2$ also coincide. It
follows that there is a cochain
$$(\phi'_\bullet,\psi'_\bullet)\in
(K^0(\gb)^1\otimes K^1(\gb)^0)\otimes I$$ such that
$$\eta=[\epsilon+ (\phi'_\bullet,\psi'_\bullet)].$$ However,
since
$(\phi'_\bullet,\psi'_\bullet)$ has coefficients in the socle
while
$\epsilon$ has coefficients in $\m_S$,
clearly $$\epsilon+
(\phi'_\bullet,\psi'_\bullet)=
\epsilon
(\phi\bullet+\phi'_\bullet,\psi_\bullet+\psi'_\bullet).$$
This proves our
assertion.\end{proof}
\begin{rem}
After \cite{atom} was written and initially posted on the arXiv,
I was informed by Michel Duflo that the 'Hard Lemma'
(Proposition 0.2) of \cite{atom}, which is at the basis of
the construction of the Jacobi-Bernoulli complex,
had been proved earlier, and used in a different Lie-theoretic context,
in the Jussieu thesis of his student Emanuela Petracci (see
Bull. Sci. Math. 127 (2003)).
\end{rem}
\newsection{Normal algebra}
Let $X\to P$ be a closed embedding of schemes.
In this section we will construct, in terms of resolutions, a
differntial graded Lie algebra $\n_{X/P}$, called the
\emph{normal algebra}, which controls motions and deformations
of $X$ in $P$.
This
will be an essential ingredient in the construction of the
abstract
tangent Lie algebra of a scheme (such as $X$). The
algebra $\n_{X/P}$ is a realization of the general \emph{normal
atom} of a subscheme (cf. \cite{atom}), one which works best if
$P$ itself is an affine or projective space. See \cite{atom}
for
more
information on Lie atoms.\subsection{Resolutions} Thus let $X$
be a
closed subscheme of a quasi-projective scheme $P$ and let
$\I=\I_{X/P}$ denote the ideal sheaf of $X$ in $P.$
Let\begin{equation}...\to F^{-1}\to F^0\to \I\end{equation} be
a
locally free resolution of $\I.$ As is well-known, $F^.$ is
locally
unique up to quasi-isomorphism. In fact, when $P$ is an
open subset of a projective space (e.g. an affine or projective
space), we may assume
\begin{equation}\label{F}F^i=N_i\O_P(-m_i):=
\bigoplus\limits_1^{N_i}\O_P(-m_i)\end{equation} for some large
integer $m_i.$   By Hilbert's
syzygy theorem the resolution $F^.$ will be finite, i.e.
bounded.
For $i=0$, each summand $\O(-m_0)$ corresponds to a generator
$f_\alpha$ of $\I$, $\alpha=1,...,N_0$. Similarly for $i<0$
each
summand $\O(-m_i)$ corresponds to a syzygy $G^i_\alpha$, that
is, a
vector
$$(g^i_{\ab})\in N_{i+1}\O(m_i-m_{i+1})$$ such that
$$d_{F^{i+1}}(\sum g^i_{\ab}G^{i+1}_\beta)=0.$$
By definition,
$$d(G^i_\alpha)=\sum\limits_\beta g^i_\ab G^{i+1}_\beta.$$
 For notational
consistency, we also set $G_\alpha^0=e_\alpha$. Thus, for each
$i\leq 0$ the $G^i_\alpha$ form a basis of $F^i.$
 Next, set $F^1=\O_P$ and let $F^._+$ be the complex in
degrees $\leq 1$ given by
$$...\to F^0\stackrel{\epsilon}{\to} F^1$$
(where $\epsilon:F^0\to F^1$ is the composite $F^0\to\I\to
\O_P$);
thus $F^._+$ is quasi-isomorphic to $\O_X[-1]$. \par
Suppose now that $P$ is open in a projective space. Note that
the
(symmetric)  multiplication map on $\O_X$ gives rise to a map
of
complexes \beq\label{gradedmult} \mu:\bigwedge\limits^2
F^._+\to
F^._+[-1]\eeq
which turns $F^.[1]$ into a \emph{differential graded
commutative $\O_P$ algebra}.
 $\mu$
is defined by descending induction, starting with the obvious
multiplication maps $$F^1\times F^1\to F^1, F^0\times F^1\to
F^0.$$
Set $G^.=\bigwedge\limits^2 F^._+$ Assuming $\mu^.$ is defined
as
map of complexes down through degree $i$, $i\leq 1$, note that
$$d^{i}_{F^._+[-1]}\circ \mu^i\circ d^{i-1}_{G}=\mu^{i+1}\circ
d^{i}_{G}\circ d^{i-1}_{G}=0$$ As $G^{i-1}$ is a free module
and
$F^._+$ is exact, this implies that $\mu^i\circ d^{i-1}_G$
factors
through a map $$\mu^{i-1}:G^{i-1}\to F^{i-2}$$ so we have
$\mu^.$ as
map of complexes down through degree $i-1$.\par For example,
for
the
map $\mu$ constructed in this way, $\mu(e_\alpha\wedge
e_\beta)$
will be the Koszul relation $f_\alpha e_\beta-f_\beta
e_\alpha$.\par Note next that,
$$\deg\mu(AB)=\deg(A)+\deg(B)-1<\min(\deg(A), \deg(B))$$
provided $\deg(A),
\deg(B)\leq 0$. Then considering the flexibility we have in
defining the various syzygy modules $F^i$, as well as $\mu$, by
descending recursion, we may  assume that\newline
\emph{(*)\ \  whenever $A,B$ are in
the standard basis that $F^.$ comes with, then so is
$\mu(AB)$;}\par
(i.e. if $AB$ yields the trivial relation, we just introduce a
basis element of the appropriate degree corresponding to the
trivial relation and assign $\mu(AB)$ to it).\par
\subsection{Definition of normal atom}
 We will denote by
$\gl^.(\I)$ the \dgla $\Hom^.(F^., F^.)$, whose term in degree
$i$
is
$$\bigoplus\limits_j\Hom(F^j, F^{i+j}),$$  whose differential
is
given by: for
$$(\phi^.)\in\gl^i(\I), \phi^j\in\Hom(F^j, F^{i+j}),$$
set $d(\phi^.)=(\psi^.),$ where
\beq\label{d-on-hom}\psi^j=d(\phi)^j=\phi^{j-1}\circ
d_{F^{j-1}}+(-1)^{i+1}d_{F^{i+j}}\circ\phi^j\eeq and whose
bracket
is
the usual graded commutator, i.e. for $(\phi^.), (\psi^.)$
homogeneous, $$[(\phi^.), (\psi^.)]=(\phi^.)\circ(\psi^.)-
(-1)^{\deg(\phi^.)\deg(\psi^.)}(\psi^.)\circ(\phi^.).$$
Note that $d\phi$ itself is the graded commutator $[(\phi^.),
d]$. There is
an
obvious action  map \beq\gl^.(\I)\times F^.\to F^.\eeq as well
as a
multiplication map induced by $\mu$: \beq F^.\times
\gl^.(\I)\to
\gl^.(\I)[-1]\eeq
 As we saw in
\cite{atom}, there is a Lie atom called the \emph{normal atom}
to
$X$ in $P$, denoted
$$\n=\n_{X/P}=\gl(\I<\O_P)=\gl(\O_P>\O_X),$$
which can be defined as the mapping cone of the evident map
\begin{equation}\label {normpair}
\gl^.(F^._+, F^1_+)\to \gl(F^1_+)\end{equation} where
$\gl^.(F^._+,
F^1_+)$
 is the subalgebra of
$\gl^.(F^._+)$ consisting of endomorphisms sending the
subcomplex
$F^1_+=\O_P$ to itself. Thus $\n$ is  also equivalent to the
mapping
cone 
\beq\label{normatom}
\gl^.(\I)\to\Hom^.(F^.,\O_P),\eeq
in other words, to the sub-dgla of $\gl^.(F^._+)$ consisting of
endomorphisms vanishing on $F^1_+.$ Simply put, we have
\begin{equation}\label{normaldgla}\n=\Hom^.(F^.,
F^._+)\end{equation}
and it is just the sub-dgla of $\gl(F^._+)$ consisting of
endomorphisms that vanish on the term (and subcomplex) $F^1_+$.
Thus
$\n$ has a dgla structure compatible with its Lie atom
structure. Note that $\n$ acts naturally on $F^._+$. Moreover
this action annihilates, hence preserves, the subcomplex
$F^1_+\subset F^._+$, hence induces an action on the quotient
$F^.=F^._+/F^1_+$. Thus we have (mutually compatible) action
pairings
$$\n\times F^._+\to F^._+, \ \ \n\otimes F^.\to F^..$$
Also, the multiplication map $\mu$ induces a multiplication
\beq
F^._+\times \n\to \n\eeq making $\n$ an $\O_X-$module.
\par In case $P$ is an affine space, because as complex (or
derived category object) $\n$ can be
identified with $\Hom^.(F^., F^._+)$ (i.e. $R\Hom(\I,\O_X)$),
it
follows that the cohomology
\begin{equation}\label{ext} H^i(\n)=\ext^i(\I,
\O_X).\end{equation}
\begin{example}\label{hypersurface}{ As a trivial example,
consider the case where
$X$ is a (Cartier) divisor in $P$, with (invertible) ideal
$\O(-X).$
Then $\gl(\I)=\O_P$ and $\n$ reduces to the complex
$\O_P\to\O_P(X)$, i.e. to $\O_X(X)[-1].$  }
\end{example}
\subsection{Relative situation} Suppose \beq\label{xyz} X\to
Y\to Z\eeq are closed embeddings. Then we can choose a
resolution $F^.$ for $\I_{X/Z}$ so that $$F^i=F^i_1\oplus
F^i_2,$$ where $F^._1$ is a resolution of $\I_{Y/Z}$, and
$F^._2\otimes \O_Y\sim_{\mathrm{qis}}F^._2\otimes F^._{1+}$ is
a
resolution of $\I_{X/Y}$. Then we have an exact triangle
\beq\label{Nxyz-1} \Hom^.(F^._2, F^._+)\to \Hom^.(F^.,
F^._+)\to
\Hom(F^._1, F^._+)\to\eeq and note that via
$$F^._{2+}\otimes\O_Y\qis F^._1\otimes F^._{1+}\qis F^._+,$$
the
first term in \ref{Nxyz-1} is quasi-isomorphic to
$$\Hom^._Y(F^._2\otimes\O_Y, F^._{2+}\otimes\O_Y)=\n_{X/Y},$$
while the third term is similarly
 quasi-isomorphic to $$\Hom^.(F^._1, F^._{1+})\otimes\O_X.$$ In
that sense then, we have a 'normal atom sequence'
\beq\label{Nxyz-2}
\n_{X/Y}\to\n_{X/Z}\to\n_{Y/Z}\otimes\O_X\to\eeq

\subsection{Jacobi complex}
As discussed in \cite{atom}, there is a \emph{Jacobi complex}
associated to $\n$, denoted $J(\n)$, which, as sheaf, lives on
$P\langle\infty\rangle$, the space of finite subsets of $P$
in
and takes the form
$$J(\n):\ \  ...\to \lambda^2\n \to\n$$
 where the term in degree $-m$, $\lambda^m\n$, is the $m$-th
 external alternating power of $\n$ and the differential is
induced
 by the bracket pairing on $\n$.
 In fact, since
$\n$ is exact off $X$, $J(\n)$ is quasi-isomorphic to, and may
be
replaced by, its (topological) restriction on $X\<\infty\>$.
The
$m$th truncation
$$J_m(\n): \ \
\lambda^m\n\to ...\to \lambda^2\n \to\n$$ lives on $X\langle
m\rangle$. The Jacobi complex $J(\n)$ admits a comultiplicative
structure
$$J(\n)\to\Sym^2(J(\n)),$$ which induces a $\C$-algebra
structure on $\C\oplus\H^0(J(\n))^*$ and on
$$R_m(\n)=\C\oplus\H^0(J_m(\n))^*.$$
Then
$$\what{R}(\n)=\lim\limits_{\leftarrow}R_m(\n)$$ coincides with
the
formal completion at $[X]$ of the local ring of Hilbert scheme
of
subschemes of $P$.\par As an alternative to working on spaces
such
as $P\langle\infty\rangle$, one may replace $\n$ by a
representative
complex of groups or $\C$-vector spaces whose terms are
acyclic,
and
in this case $J(\n)$ can be replaced by the standard Jacobi
complex
where $\lambda^i$ is just the ordinary exterior power over
$\C$.
Such a representative can be taken to be the sections of an
injective or soft resolution. However in most applications
considered here, $X$ is an affine scheme and it suffices to
replace
$\n$ simply by the corresponding complex of global sections. We
will
refer to this as the \emph{non-sheafy} version of $J(\n)$ in
the
affine case.\par
\subsection{Regular embeddings} One point which has to be noted
is that if
$X\to P$ is a regular embedding (of codimension $>1$), our
normal
atom differs from the usual normal sheaf $$N=\Hom(\I, \O_X).$$
In
fact, since $\n$ is  $\RHom^.(\I,\O_X)$, it is given
$\Hom^.(F^.,\O_X)$, where $F^.$ is the Koszul complex. Since
the
differentials of $F^.$ have coefficients in $\I$, it follows
that
$\n$ is the complex with trivial differentials and terms
\begin{eqnarray*}\n^i={\bigwedge\limits^i}_{\O_X}N,
i=1,...,\codim(X,P)\\=0 \ \ \ {\rm otherwise.}\end{eqnarray*}
Why is
it sufficient to work just with $N$? We may reason as follows.
The
bracket of $\n$, restricted on $N$ (which itself is viewed as a
complex, indeed a subcomplex of $\n$, via an injective
resolution)
has 2 components: one is a map \beq\label{Nbracket}\Sym^2N\to
N[-1]\eeq while the other is a (locally determined) map
$$\Sym^2(N)\to \n^2=\bigwedge\limits^2 N$$ (note that
$\bigwedge\limits^2$ becomes $\Sym^2$ on $\n^1$). We claim the
second map vanishes. Indeed, locally, a section of of
$\bigwedge\limits^iN$ lifts to a section of $F^{-(i-1)}\check{\
\
}$, and this yields an endomorphism of degree $i$ of $F^._+$
via
interior multiplication. Moreover, (signed) commutator of
endomorphisms, which yields the bracket on $\n$, corresponds to
wedge product. Thus for local sections $a,b$ of $N$, the
commutator
is $ab+ba$ which corresponds to $a\wedge b+b\wedge a=0$. This
corresponds to the well-known fact that 'local obstructions to
deforming complete intersections vanish'. Therefore $N$ with
induced
bracket (\ref{Nbracket})  is a sub-Lie atom of $\n$. This
remark
will be generalized below via the notion of \emph{reduced}
normal
atom.
\subsection{Weak equivalence}The following definition will be
useful
\begin{defn}A homomorphism
\beq\label{weak}\theta:\g_1\to\g_2\eeq of
Lie atoms or dg Lie algebras is said to be a \emph{direct weak
equivalence} or \emph{diweq} if it induces an
isomorphism on $H^i, \forall i\leq 1$,
and an injection on $H^2.$ A general \emph{weak equivalence} or
\emph{weq} is a
composition of direct weak equivalences and their inverses. The
relation of \emph{weak equivalence} is the equivalence relation
generated by the existence of a direct weak equivalence.
\end{defn} The discussion of the preceding
subsection proves \begin{lem} If $X\to P$
is a regular embedding then the normal atom $\n_{X/P}$ is
weakly
equivalent to the shifted normal bundle $N_{X/P}[-1]$ with
trivial
bracket.\end{lem} The following observation is easy to prove
but
will be quite useful in what follows.
\begin{lem} A weak equivalence (\ref{weak}) induces an
isomorphism on
deformation rings $$R_m(\g_1)\simeq R_m(\g_2).$$\end{lem} For
the
normal atom $\n$ of a regular embedding, clearly $N[-1]\to\n$
is
a
weak equivalence.\par
 \begin{rem}{\rm If $X$ is smooth and $P$
 is an open subscheme of a projective space, it is possible to
 identify the Lie atom $N$ as constructed above with the
 'differential'
normal atom  $N^{\rm diff}$, defined by the complex
$$T_{P/X}\to T_P\sim N^{\rm
diff}[-1],$$ where $T_{P/X}$ is the algebra of vector fields on
$P$
tangent to $X$ along $X$ (see \cite{atom}). Briefly, the
methods
of
the next section show that we can define a diagram
\begin{equation}\begin{matrix}
T_{P/X}&\to &T_P\\ \downarrow&&\downarrow\\
\gl(\I)&\to&\Hom^.(F^.,\O_P)\end{matrix}\end{equation} by
making
vector fields act on the equations of $X$. This diagram
induces
a
quasi-isomorphism of Lie atoms $$N^{\rm diff}[-1]\to N.$$ In
this
sense then, $\n$ is weakly equivalent to $N^{\rm diff}[-1].$
}\end{rem} \subsection{Reduced normal atom}\label{rednormalsec}
Consider now the case where $P$ is an open subscheme of an
projective space, but $X\to P$ is arbitrary.
Our purpose here is to show that the normal atom
$\n$
can be replaced by another, weakly equivalent dgla $\n_\red$
which is 'smaller' and
more
convenient for many applications. Essentially, $\n_\red$ is
obtained from $\n$ by removing the trivial Koszul syzygies,
which play no role in deformation and obstruction theory.
This construction is motivated by, and analogous to, the
\emph{cotangent complex} of a scheme, which goes back to
Grothendieck.
\par
Consider a resolution $F^.\stackrel{\epsilon}{\to}I$ as above,
where $F^0$ is a free module on a basis $(e_\alpha)$ that maps
to a system of generators $(f_\alpha)=(\epsilon(e_\alpha))$ of
$I$.
Let $K_0\subset F^0$ denote the 'Koszul' submodule, generated by
all elements of the form
$$e_{\alpha\beta}=f_\alpha e_\beta-f_\beta e_\alpha.$$ Note that
if the $e_\alpha$ form a regular sequence, then
$K_0=\ker(\epsilon)$.  Let
\beq\label{K_0}...\to F^{-2}_0\to F^{-1}_0\to K_0\eeq be a free
resolution of $K_0$, where $F_0\inv$ is a free module on
generators $s_{\alpha\beta}$ that map to  $e_{\alpha\beta}$. We
may assume
our resolution $F^.$ of $I$ is constructed in such a way that
for each $i<0$, there is a splitting of $F^i$ into free modules
$F^i=F^i_0\oplus F^i_1$. Such resolutions $F^.$ will be called
\emph{standard}. Thus any resolution is equivalent to a standard
one. Note that, for a standard resolution $F^.$, the complex
$F^{<0}$, itself a quotient of $F^.$ as well as $F^._+$, fits in
an exact sequence
\beq\label{F_0FF_1}0\to F^{<0}_0\to F^{<0}\to F^{<0}_1\to 0.\eeq
If the $e_\alpha$ form a regular sequence, then clearly
$F_1^{<0}=0.$\par
Now the idea is that if we deform the subscheme $X\subset P$ by
deforming its resolution, one might as well start with a
standard resolution and deform the latter in the more restricted
universe of standard resolutions, rater than all resolutions.
This is accomplished by the dgla $\n_\red$ that we now
introduce.\par
 Start with a standard resolution $F^.$ as above.
Then define a map \beq\label{pidef1}\pi:\Hom(F^0, F^i_+)\to
\Hom(F\inv_0, F^{i-1}_+)\subset \Hom(F\inv, F^{i-1}_+), \forall
i\leq 1\eeq by the formula
\beq\label{pidef2}\pi(e_\alpha^*\otimes g)
=(-1)^i\sum\limits_\beta
s_{\alpha\beta}^*\otimes\mu(e_\beta\otimes g),\forall g\in
F^i\eeq where $$\mu:F^0\otimes F^i_+\to F^{i-1}_+$$ is the
graded multiplication (induced by the multiplication on $A_X$ as
in \refp{gradedmult}.). One can easily check that
$$s_\ab=\mu(e_\alpha\otimes e_\beta).$$ In particular,
$\mu(F^0\otimes F^0)\subset  F\inv_0.$
Now the commutative diagram (with left vertical maps
$\id\otimes d$)
$$\begin{matrix} F^0\otimes F_+\inv&&\\ \downarrow&&\\
F^0\otimes F^0_+&\to&F\inv\\
\downarrow&&\downarrow\\ F^0\otimes F^0_+&=&F^0\end{matrix}$$
shows that the composite
$$F^0\otimes F\inv \stackrel{\id\otimes d}{\to} F^0\otimes
F^0\stackrel{\mu}{\to} F\inv_0$$ factors through the kernel of
$F\inv_0\to K_0$, hence
we may assume that
$\mu|_{F^0\otimes F\inv}$ factors through $F^{-2}_0$.
Then again the commutative diagram
$$\begin{matrix} F^0\otimes F^{-2}&&\\\downarrow&&\\
F^0\otimes F\inv&\to&F^{-2}_0\\
\downarrow&&\downarrow\\ F^0\otimes
F^0&\to&F\inv_0\end{matrix}$$
shows may assume that
$\mu|_{F^0\otimes F^{-2}}$ factors through $F^{-3}_0$, etc.
Thus,
 we may assume
\beq\label{mutoF0}\mu(F^0\otimes F^i)\subset F^{i-1}_0,\l
\forall i\leq 0,\eeq and consequently
\beq\label{pitoF0}\pi(\Hom(F^0, F^i_+))\subset \Hom(F\inv_0,
F^{i-1}_0)\eeq
We consider $\pi$ as an endomorphism of $\n$. Note that,
trivially, $\pi^2=0$.
\begin{lem}\label{pi-homo} $\pi$ is a Lie
homomorphism.\end{lem}\begin{proof}
It suffices to check this on $\Hom(F^0, F^i_+)$ and the only
nontrivial case is $i=0$, in which case we must prove
$$\pi([e_{\alpha_1}^*e_{\gamma_1},
e_{\alpha_2}^*e_{\gamma_2}])=[\pi (e_{\alpha_1}^*e_{\gamma_1}),
\pi (e_{\alpha_2}^*e_{\gamma_2})].$$
Suppose first that $\gamma_1\neq\alpha_2, \gamma_2\neq\alpha_1$,
so the LHS is clearly zero.
Then we compute the RHS:
$$\sum\limits_{\beta_1,
\beta_2}[s_{\alpha_1\beta_1}^*s_{\beta_1\gamma_1},
s_{\alpha_2\beta_2}^*s_{\beta_2\gamma_2}]=
s_{\alpha_1\alpha_2}^*s_{\gamma_1\gamma_2}-
s_{\alpha_2\alpha_1}^*s_{\gamma_2\gamma_1}=0$$
If $\gamma_1=\alpha_2, \gamma_2\neq\alpha_1$ then the RHS is
$$\sum\limits_\beta s_{\alpha_1\beta}^*s_{\beta\gamma_2}$$ which
clearly equals the LHS. Other cases are similar.
\end{proof}
\begin{lem}\label{piproj} The map $$i_{F\inv_0}^*\circ
d:\Hom(F^0, F^i_+)\to \Hom(F\inv_0, F^i)$$ coincides with
$$i_{F\inv_0}^*\circ(-d\circ\pi+\pi\circ d).$$\end{lem}
\begin{proof} Pick $\phi=e_\alpha^*\otimes g\in\Hom(F^0, F^i)$.
Then we have
$$d\phi=\phi\circ d-(-1)^id\circ\phi=-s_\ab^*f_\beta\otimes
g-(-1)^i e_\alpha^*\otimes dg.$$
Hence, the component of $\pi d\phi$ in $\Hom( F\inv, F^{i})$ is
$$(\pi d\phi)^{1,i+1}=s_\ab^*\otimes \mu(e_\beta\otimes dg).$$
On the other hand, since
$$\pi(\phi)=(-1)^is_\ab^*\mu(e_\beta\otimes g),$$
we have
$$d\pi(\phi)^{1,i}=s_\ab^*\otimes(f_\beta g+\mu(e_\beta\otimes
dg))$$ Therefore,
$$d\phi (s_\ab)=-d\pi\phi(s_\ab)+\pi d\phi(s_\ab), \forall
\alpha, \beta$$ which proves our assertion.
\end{proof}
We exploit the lemma as follows. Set $$d'=d+d\pi-\pi d
:\n\to\n.$$ Then it is immediate that $$(d')^2=d\pi d-d\pi
d=0.$$ Therefore we may define a complex
$$\n'=(\n, d').$$ By the very definition of $d'$, the map
$$\id+\pi:\n'\to\n$$ is a map of complexes. The map is clearly
an isomorphism, with inverse $\id-\pi$. Now define subgroups of
$\n'$ as follows:
\beq\label{n_red}\nred =\bigoplus\limits_{i>0,j}\Hom(F_1^{-i},
F^j_+)\oplus\bigoplus\Hom(F^0, F^j_+),\eeq
\beq\label{n_0}\n_0=\bigoplus\limits_{i>0, j}\Hom(F^{-i}_0,
F^j_+).\eeq
Note that when the $e_\alpha$ form a regular sequence,
$$\nred =\Hom^.(F^0, F^._+)\sim \Hom(F^0, A_X)[-1].$$
Now  Lemma \ref{piproj} means precisely that the obvious
splitting of $\n$ as group yields a splitting of the
\emph{complex} $\n'$:
\beq\label{n'split}\n'=\nred\oplus \n_0.\eeq Now, it is easy to
check that $\nred$ is in fact a sub-dgla of $\n'$. I claim next
that $\id+\pi$ is a Lie homomorphism on $\nred$, i.e. that
\beq\label{id+pi}[\phi_1+\pi(\phi_1),\phi_2+\pi(\phi_2)]=
[\phi_1,\phi_2]+\pi([\phi_1, \phi_2]).\eeq
Indeed it suffices to check this when $\phi_1, \phi_2$ are
bihomogeneous, in which case,
by \refp{pitoF0}., it holds for trivial reasons except when
$\phi_1,\phi_2\in\Hom(F^0, F^0)$. In this case, \refp{id+pi}.
follows from the fact that $\pi$ is a Lie homomorphism (Lemma
\ref{pi-homo}). The same lemma also shows that $d'$ is a
derivation, so that $\n'$ is a dgla.\par Now, since
$$\n_0\sim\Hom^.(F^{<0}_0, A_X[-1]),$$ this complex is acyclic
in degrees
$<2$. Therefore the inclusion $\n_\red\to\n'$, hence also the
induced map $\nred\to \n$, are weak equivalences.\par
Also, it is easy to check that the usual argument showing that
the resolution $F^.$ is uniquely determined up to homotopy also
shows that $\nred$ is uniquely determined up to dgla homotopy.
\par We summarize the foregoing discussion in the following
statement. \begin{thm} There exists a dgla $\nred=\nred_{X/P}$,
canonical up to dgla homotopy, together with a weq $$\nred\to
\n_{X/P},$$ such that, whenever $X\to P$ is a regular embedding,
we have $$\nred\sim N_{X/P}=\Hom(I_X, A_X)[-1].$$
\end{thm}
Note that if a subset of the generators $f_\alpha$ of $I$ form a
regular sequence, they contribute to $\n_\red$ only via
$\Hom(F^0, F^._+)\sim\Hom(F^0, A_X)$. Combining
this with \refp{Nxyz-2}. above, we conclude the
following
\begin{cor}\label{nred-xyz} If $X\to Y\to Z$ are closed
embeddings with $Y\to Z$
regular, then there is an exact sequence $$ 0\to\nred_{X/Y}\to
\nred_{X/Z}\to N_{Y/Z}[-1]\otimes A_X\to 0.$$\end{cor} Next,
note that if
$f\in A_P$ and $\bar{f}\in A_X$ is its image,
then we have a closed embedding of principal open affines
$$X_{\bar{f}}\to P_f$$ and a resolution of $X\to P$ yields one
of $X_{\bar{f}}\to P_f$, therefore we have a dgla
quasi-isomorphism
$$(\nred_{X/P})_{\bar{f}}\sim \nred_{X_{\bar{f}}\to P_f}.$$ In
particular, if $X_{\bar{f}}\to P_f$ is regular, then
$\nred_{X/P}$ is acyclic in degree $>1$ over $X_{\bar{f}}.$ So
we conclude \begin{cor} For $i>1,\  H^i(\nred)$ is supported on
the set of points where $X\to P$ is not a regular
embedding.\end{cor}

\begin{rem} To fix ideas, we have worked in this section in the
setting of coherent sheaves. However when $P$ is an affine or
projective space, there is no difficulty in shifting from the
sheafy setting to the module (resp. graded module) setting,
using the standard exact functors $\Gamma$ (in the affine case)
or $\Gamma_*$ (in the projective case).\end{rem}
\newsection{Affine schemes: Tangent
atom}\label{affine}\setcounter{equation}{0} Our
purpose eventually is to define, for any algebraic scheme $X$,
a
Lie
atom $\t_X$ called the \emph{tangent atom} of $X$, that will
play a
role in deformation theory akin to that of tangent algebra for
a
smooth scheme.
 We begin in this section with the case where $X$ is specified
as
embedded as closed subscheme in a fixed
affine space $P=\A^n$.
The complex $\t_X=\t_X(P)$ that we construct will depend on the
embedding $X\to P$- as well as a number of other choices- but,
as we
shall see, the weak equivalence class of $\t_X(P)$ will be
independent of all the choices, and will depend--
functorially--
 only on
the isomorphism class of $X$.
Moreover, $\t_X(P)$ will in fact be
(weakly equivalent to) a dgla, rather than a general Lie
atom.\par
To fix ideas, we will be working in this section in the setting
of (finitely generated) modules over the affine coordinate ring
$A_P$. There is no difficulty in carrying the construction to
the
setting of coherent sheaves over $P$, but we will not need
this.
In any event, the extension to the global (non-affine) setting
is far from straightforward and will be taken up in the next
section.
 \subsection{Construction of differential graded group} In this
subsection we will construct the differential graded group
underlying the tangent dgla $\t_X$. The proof that this forms a
complex, as well as the construction and properties of the Lie
bracket, will be deferred to the next subsection.\par
Now our construction of the tangent atom $\t_{X}$ is motivated
by
the fact that for $X$ smooth, its tangent algebra is
quasi-isomorphic to the complex (of Lie atoms)
\begin{equation}\label{tangcx}\I_X\otimes T_P\to T_P\to
N_{X/P}.
\end{equation} There is another Lie algebra denoted $T_{X/P}$
of
ambient vector fields on $P$ tangent to $X$ along $X$ (or what
is the same, derivations on $A_P$ taking the ideal $I_X$ to
itself), which itself is the mapping cone of $T_P\to N_{X/P}$,
so that $T_X$ is the mapping cone of $I_X\otimes
T_P\to T_{X/P}$.
The plan is to construct a complex like
(\ref{tangcx}) in the general case, with $\I_X$ replaced,
naturally
enough, by a free resolution of itself, and $N_{X/P}$ replaced
by (a
convenient model of) the shifted normal atom $\n_{X/P}[1]$, or
its reduced version $\nred_{X/P}[1]$.
Indeed,
everything but the second arrow has already been constructed.
Thus
it remains to construct (something like) the second arrow,
though in
the singular case one cannot expect this arrow to be
surjective.
\par Thus, the tangent atom
$$\t=\t_X=\t_X(P)$$ will be a complex $\t^.$, depending on the
embedding $X\to P$. Alongside $\t$, we will be constructing
a complex $\t^._{X/P}$ analogous to $T_{X/P}$ in the smooth
case.
 To define $\t$, fix a (finite) free resolution
$F^.$ in degrees $\leq 0$ of the ideal $I=I_{X/P}$ (non-sheafy
version), as above, with
$F^._+$ the corresponding resolution in degrees $\leq 1$ of
$A_X$,
with $F^1_+=A_P$. Recall that the normal atom
$$\n=\n_{X/P}=Hom^.(F^.,F^._+).$$ Thus
$$\n^i=\Hom^i(F^.,F^.)\oplus
\Hom(F^{-i+1},A_P)=\gl^i(I)\oplus\check{F}^{-i+1},$$
where, as before, we identify
\beq\gl^i(\I)=\gl^i(F^.)=\bigoplus\limits_j\Hom(F^j,
F^{i+j}).\eeq We view $\n$ as the sub-dgla of $\gl^.(F^._+)$
consisting of elements annihilating the subcomplex $F^1_+$ of
$F^._+$.
 Now the terms of $\t^.$ are defined as
follows $$\t^i=T_P\otimes F^{i+1}_+\oplus\n^i.$$ More
explicitly,
\begin{eqnarray}
\t^i&=&(T_P\otimes F^{i+1}_+)
\oplus \gl^i(I),\ \  i\leq 0,\\
&=&\gl^i(I)\oplus\Hom(F^{-i+1},A_P),\ \ \ \ \   i\geq 1.
\end{eqnarray}
So, schematically,\beq
\t_X:...\to\begin{matrix}\gl^{-1}(I)&&\gl^0(I)&&\gl^1(I)\\
\oplus&\to&\oplus&\to&\oplus\\
(T_P\otimes F^0)&&T_P&&\check{F}^0\end{matrix} \to ...
\eeq (with the middle term in degree 0). Thus as graded group,
$\t_X$ in the affine case is just $\n\oplus (T_P\otimes
F^._+)$.
 The analogous complex
$\t^._{X/P}$ is constructed similarly, with $T_P\otimes F^._+$
replaced by $T_P$, i.e.$$\begin{matrix}\t_{X/P}^i=&\n^i,\ \ \
i\neq
0\\ &\n^i\oplus T_P, i=0.\end{matrix}$$
.\par
Note that all $\t^i$ are finitely generated locally free (even
globally free, if $P$ is affine space) sheaves, and all but
finitely
many vanish. We will define the differentials on $\t^.$ in such
a
way that all components with source (resp.
target) in $\n$ (resp. $T_P\otimes F^._+$) will coincide with
the
original ones from $\n$ (resp. $T_P\otimes F^._+$), so that
there
will be a exact sequences of complexes $$0\to\n\to\t\to
T_P\otimes
F^._+\to 0,$$
$$0\to\n\to\t_{P/X}\to T_P\to 0.$$
Thus what really needs to be defined is the 'mixing'
part, i.e. the differential components going from $T_P\otimes
F^._+$
to $\n$.\par
 Now explicitly, the various components
 of the differentials $d^i:\t^i\to\t^{i+1}$ are defined as
follows. \begin{itemize}
\item For $i=0$, the components of $d^0$
are:\begin{itemize}\item the map
$\gl^0( I)\to\gl^1( I)$ is the differential of $\gl^.(
I)$;\item
the
map $\gl^0( I)\to \Hom(F^0,  A_P)$ sends a graded map
$(\phi^.)$
to
$\epsilon\circ\phi^0$;\item let $(e_\alpha)$ be a basis
for the free module $F^0$ and $f_\alpha=\epsilon(e_\alpha)$
the corresponding generator of $I$, so that
$$\epsilon=\sum_\alpha f_\alpha e_\alpha^*;$$
then the map $$T_P\to \Hom(F^0,
 A_P)$$ sends a
section (vector field) $v$ to
\begin{equation}d^0(v)=(v(f_a))=\sum_\alpha
v(f_\alpha)e^*_\alpha\end{equation} \item the map
$$u=(u^.):T_P\to\gl^1( I)= \bigoplus\limits_j\Hom(F^{j}\to
F^{j+1})$$ has $j$-th component $u^j$ defined as follows. Let
$(G_\alpha^{j+1})$ be a basis of $F^{j+1}$ and $(G_\beta^j)$ be
a basis of
$F^j$. Thus each $G_\beta^j$ is a syzygy between the
$G_\alpha^{j+1}$ and
may be written as a vector $G_\beta^j=(g^j_{\beta\alpha})$ with
each $g^j_{\beta\alpha}\in A_P$ and
\begin{eqnarray}\label{syz}\sum\limits_\alpha
g_{\beta\alpha}^jd^{j+1}(G^{j+1}_\alpha)=
\sum\limits_\alpha
g_{\beta\alpha}^j g^{j+1}_{\alpha\eta} G^{j+2}_\eta=
0,& j+1<0\\
\sum\limits_\alpha g_{\beta\alpha}\epsilon(e_\alpha)=
\sum\limits_\alpha g_{\beta\alpha}f_\alpha=0,& j+1=0
\end{eqnarray} In fact the second equation becomes a special
case of the first, for $j+1=0$ if we take for $(G^1_\eta)$ the
standard basis $(1)$ of $A_P=F^1_+$. Note that \refp{syz}. is
equivalent to the scalar identity
 \begin{equation}\label{syzy} \sum\limits_\alpha
g_{\beta\alpha}^j g^{j+1}_{\alpha\eta}=
0, \forall \beta, \eta.\end{equation}
Then set
$$u^j(v)=\sum\limits_{\alpha,
\beta}v(g^j_{\beta\alpha})G^{j*}_\beta
G^{j+1}_\alpha.$$ Thus, all in all, the differential
$$d:T_P\to\n^1$$ is defined by
\beq\label{d-on-TP}d(v)=\sum\limits_{j\leq 0, \alpha,
\beta}v(g^j_{\beta\alpha})G^{j*}_\beta
G^{j+1}_\alpha.\eeq This may be explained by noting that $v\in
T_P$ defines in an obvious way a linear map $I_X\to
A_X=A_P/I/X$. Then $d(v)$ is just the extension of this to a
map
of resolutions $F^.\to F^._+$. Another explanation is this.
Since $F_+^.$ is endowed with a specified basis, the derivation
$v$ of $A_P$ extends to an additive termwise endomorphism (in
fact, an $A_P$-derivation compatible with the derivation $v$ of
$A_P$) of $F_+^.$; then $d(v)$ is just the commutator of $v$
with the differential $d$, i.e.
\beq\label{dTXP-def-bracket} d(v)=[v,d]\eeq
(note that a priori $[v,d]$ vanishes on $F^1_+$ because $d$
does, $v$ kills 1, and $[v,d]$ is linear, as the commutator of a
derivation and a linear map).
\end{itemize}\item for $i<0$, the components of $d^i$
 are:\begin{itemize}
\item the map $\gl^i( I)\to\gl^{i+1}( I)$ is the
$\gl^.( I)$-differential;\item the map $F^{i+1}_+\otimes T_P\to
F^{i+2}_+ \otimes T_P$ is $d_{F^._+}\otimes\id_{T_P}$;\item the
map
$$u=(u^.):F^{i+1}\otimes
T_P\to\gl^{i+1}( I)=\bigoplus\limits_j\Hom(F^j, F^{j+i+1})$$ is
given
by \begin{equation}\label{J}u^j(G\otimes
v)=\sum\limits_{\alpha,
\beta}
v(g^j_{\beta\alpha})G^{j*}_\beta\mu(G\otimes
G^{j+1}_\alpha)\end{equation}
where $$\mu:
F^a_+\otimes F^b_+\to F^{a+b-1}_+$$ is the pairing obtained as
the extension of the multiplication map on $A_X$ to its free
$A_P$-resolution (cf. \refp{gradedmult}.).\par Note that,
compared to the definition of $dv\in\n^1$ above, we can write
this component in the form $\mu(G.dv)$ where the $\mu$
multiplication refers to the natural extension of $\mu$
$$\mu:F^k\times\n^i\to\n^{i+k-1}$$$$G^k_\gamma\times
G^{j*}_\beta\otimes G_\alpha^{j+i}\mapsto
G^{j*}_\beta\otimes\mu(G^k_\gamma G^{i+j}_\alpha).$$
As in \refp{dTXP-def-bracket}., this definition may be
explained
by the remark that $G\otimes v$ defines a derivation of degree
$i-1$ of $F^._+$ as $A_P$ module, and $d(G\otimes v)$ is just
the commutator:
\beq \label{dTX-def-bracket}d(G\otimes v)=[G\otimes v, d]\eeq
\end{itemize}
\item For $i>0$:\begin{itemize} \item The map
$\gl^i( I)\to\gl^{i+1}( I)$ is the $\gl^.(
I)$-differential;\item
the map $\gl^i( I)\to\Hom(F^{-i}, A_P)$ sends a graded map
$(\phi^.)$ to $\epsilon\circ\phi^i$ where $\phi^i:F^{-i}\to
F^0$
is
the degree-$i$ component of $\phi$.\end{itemize} All other
components are zero.\end{itemize}This completes
the construction of the differential on $\t_X$; we will shortly
prove that the square of this differential is zero, making
$\t_X$ a complex.
\par The differentials of the
complex $\t_{X/P}$ are defined
compatibly with $\t_X$, so that $\t_{X/P}$ will be a subcomplex
of $\t_X$.
Thus (modulo verifying that $\t_X, \t_{X/P}$ are complexes- see
below),
 there is a short exact sequence of
complexes (exact triangle)
 \beq\label{ITP->TP/X} 0\to T_P\otimes F^.\to \t_{X/P}\to
\t_X\to 0\eeq
\begin{example}\label{quadcone}{
Continuing with example \ref{hypersurface}, suppose $X$ is a
hypersurface with equation $f$ of degree $d$ in $P=\C^n$. As we
have
seen, $\n_{X/P}$ can be identified cohomologically  with $$ A_P(d)/f A_P\simeq
A_X,$$ or more precisely with the complex
$$A_P\stackrel{f}{\to} A_P(d).$$ The tangent algebra $\t_X$ is the complex
$$T_P(-d)\stackrel{(L_f, -e_f)}{\to} T_P\oplus A_P
\stackrel{(e_f, L_f)}{ \to} A_P(d)$$
where $L_f$ is left (post) multiplication by $f$ and $e_f$ is evaluation on $f$.
Thus the map\nl $T_P\to A_P(d)$ sends a vector field
$\partial/\partial
X_i$ to $\partial f/\partial X_i$. \par If $f$ is nondegenerate,
this
map
yields an injection $$T_P/ fT_P\to \n_{X/P}=A_P(d)/fA_P.$$ Thus,
in this case, $\t_X$ can be
identified cohomologically with a shift of the Milnor algebra
$$A(f)=\C[X_1,...,X_n]/(f,\partial f/\partial X_1,...,\partial
f/\partial
X_n)$$
 (which is
finite-dimensional iff $X$ has isolated singularities).}
\end{example}

\subsection{Construction and basic properties of dgla} In this
subsection we will construct $\t_X$ and $\t_{X/P}$ as dgla's
and
establish their basic properties and interrelation.\par  In
what
follows we shall be working in
the derived category of coherent sheaves on $P$, and in
particular identify a
sheaf of modules, e.g. $ A_X$, with a free resolution of
itself.
Also, if no confusion is likely, we shall usually make no
distinction between a finitely generated module and the
associated coherent sheaf, e.g. $A_X$ and $ A_X$.\par The main
results of this subsection are summarized in the following.
\begin{thm}\label{affine tx}\begin{enumerate}
\item $\t^._X$ and $\t_{X/P}^.$ as defined above are complexes
and admit structures of dgla; as such, both act on $ A_X$;
$\t_X$ is also an $ A_X$- module.\item $\t_{X/P}$ is a
subcomplex of $\t_X$ and The natural map $\t_{X/P}\to\t_X$ is
a dgla homomorphism whose mapping cone is quasi-isomorphic as
Lie atom  to $ I[1]\otimes T_P=F^.[1]\otimes T_P$.
\item Endowing $F^.\otimes T_P\sim I\otimes T_P$ with its
natural Lie algebra structure as subalgebra of $T_P$, there is
a
Lie homomorphism $J: F^.\otimes T_P\to \t_{X/P}$ which realizes
$F^.\otimes T_P$ as Lie ideal in $\t_{X/P}$, such that the
quotient (= mapping cone of $J$) is $\t_X[-1]$.
\item
The weak equivalence class of $\t^._X$ is independent of the
embedding\nl $X\to P$ and is functorial in $X$, in the
following
sense. To a morphism of affine schemes $$f:X\to Y$$ is
associated a map of Lie atoms $$df:\t_X\to f^!\t_Y$$ such that
the pair $$X\mapsto {\rm{ weak\  equivalence\  class\  of\ }}
\t_X$$
$$f\mapsto df$$ is a functor from the category of affine
schemes
over $\C$ to that of weak equivalence classes of dgla's and
dgla
homomorphisms over $\C$.
\end{enumerate}\end{thm}
\begin{proof} The proof is somewhat long, so we break it into
steps.
\subsubsection{$\t^._{X/P}$ is a complex}
\  \nl This means
 $d^{i+1}\circ d^i=0.$ For all $i\neq 0$ this vanishing
  is obvious from the
analogous property for the normal atom $\n^.=\n^._{X/P}$. For
$i=0$, the vanishing of $d^1\circ d^0$ on the $\gl^0=\n^0$
 summand of $\t_{X/P}^0$
 similarly follows from the
corresponding fact for the normal atom $\n^.$.\par
For the other summand $T_P$ of $\t_{X/P}^0$,
consider a section $v$. Then $d^1\circ d^0(v)$ has
components in
$\Hom(F^i, F^{i+2}_+)$ for all $i\leq -1$, and we claim all of
them vanish.
We focus first on the case $i=-1$, i.e. we will show that the
component of
$d^1\circ d^0(v)$ in $\Hom(F^{-1}, F^1_+)=\Hom(F^{-1},  A_P)$
vanishes, as the proof for
the other components is similar. The latter component is a sum
of 2
terms, corresponding to 2 components $d^{00}(v), d^{01}(v)$ of
$d^0(v).$ The first component $d^{00}(v)$ is the map $$F^0\to
A_P,$$
$$e_\alpha\mapsto v(f_\alpha)$$
 taking the basis element $e_\alpha$ corresponding to the
generator
$f_\alpha$ to $v(f_\alpha)$. Then the component of
$d^1(d^{00}(v))$
in $\Hom(F^{-1}, A_P)$ is the map $$F^{-1}\to  A_P$$
$$G^{-1}_\beta=(g^{-1}_{\beta\alpha}) \mapsto
\sum\limits_\alpha
g^{-1}_{\beta\alpha}v(f_\beta)$$ where
$G^{-1}_\beta=(g^{-1}_{\beta\alpha})$
 is a basis
element of $F^{-1}$ which, being a syzygy, satisfies
$$\sum\limits_\beta g\inv_{\ba}f_\beta=0.$$
 The other component $d^{01}(v)$ of $d^0(v)$ is the
map
$$F^{-1}\to F^0,$$
$$G\inv_\beta \mapsto \sum\limits_\alpha
v(g\inv_{\ba})e_\beta$$
and then the component of $d^1(d^{01}(v))$ in $\Hom(F^{-1},
A_P)$ is
the map
$$G\inv_\beta\mapsto \sum\limits_\alpha
v(g\inv_{\ba})f_\beta.$$
 Thus in total, the component
of $d^1d^0(v)$ in $\Hom(F^{-1}, F^0)$ is the map
\begin{eqnarray*}G\inv_\beta\mapsto\sum\limits_\alpha
g\inv_{\ba}v(f_\beta)+
\sum\limits_\alpha v(g\inv_{\ba})f_\beta\\
=v(\sum\limits_\alpha g\inv_{\ba}f_\beta)=0.\end{eqnarray*} The
vanishing of the component of $d^1\circ d^0(v)$ in
$\Hom(F^i, F^{i+2}_+)$ for each $i<-1$ follows
similarly from the derivation property of $v$
and the characteristic property of syzygies \refp{syzy}.:
thus, the component in question is just
\begin{eqnarray*}\sum\limits_{\beta,\alpha,\eta}(g^i_\ba
v(g^{i+1}_{\alpha\eta})+v(g^i_\ba)
g^{i+1}_{\alpha\eta})G^{i*}_\beta G^{i+2}_\eta\\
=\sum\limits_{\beta,\alpha,\eta}v(g^i_\ba
g^{i+1}_{\alpha\eta})G^{i*}_\beta G^{i+2}_\eta = 0
\end{eqnarray*}
This completes the
verification that
$\t^._{X/P}$ is a
complex.\par
A more abstract, less explicit proof that $d^2(v)=0$,  based on
\refp{dTXP-def-bracket}., can be given as follows:
$$d^2(v)=[[v,d],d]=[v,[d,d]]-[d,[v,d]]$$ the latter by the
Jacobi identity. Now $[d,d]=d^2/2=0$ and
$[d, [v,d]]=[[v,d],d]$ because $d$ and $[v,d]$ both have degree
1. Therefore $$[[v,d],d]=-[[v,d],d],$$ hence
$d^2(v)=[[v,d],d]=0.$\qed
\subsubsection{$\t_X$ is a complex}\ \par
Now given that that $\t^._{X/P}$ is a
complex, proving the same for $\t^._X$ amounts to showing that
$d^{k}\circ d^{k-1}$ vanishes on
the summand
$F^k\otimes T_P$ of $\t_X^{k-1}$,
 for all $k\leq 0$. To bring out the
 idea, we will first work out the case $k=0$ as other
 cases are only notationally more complicated. Thus for a
section $s=e_\beta\otimes v$ of $F^0\otimes T_P$,
$d^{-1}(s)$ has 2 components: one in $T_P$, viz.
\begin{center}$f_\beta v\in T_P,\ \ $ and \end{center}
$$-\bigoplus\limits_{i\leq 0,\gamma,\alpha}g^{i*}_\gamma
v(g^i_{\alpha\gamma})\mu(e_\beta g^{i+1}_\alpha)
\in\bigoplus\Hom(F^i, F^i).$$ Now applying $d^0$ to this, we
get
components in $\Hom(F^i, F_+^{i+1})$ for all\nl $i\leq 0$, and
we claim
they are all zero There are 2 cases. First,
$i=0$. This component gets just 2 contributions, one from
$T_P$,
viz
$$\sum\limits_\gamma e_\gamma^*f_\beta v(f_\gamma),$$ and one
from
$\Hom(F^0,F^0)$ (via post-composing with $d_{F^._+}$), which is
just
$$-\sum\limits_\gamma e_\gamma^*f_\beta v(f_\gamma),$$ so the
total is zero.\par Next, take $i<0$. Then we get 3
contributions
in
$\Hom(F^i, F^{i+1}))$: one from $T_P$, which equals
$$I:
\sum\limits_{\gamma,\alpha} f_\beta
v(g^i_{\gamma\alpha})G^{i*}_\gamma
G^{i+1}_\alpha, $$ one from $\Hom(F^i,F^i)$
(via post-composing with $d_{F^.}$), which equals
$$II:-\sum\limits_{\gamma,\alpha}G^{i*}_\gamma
v(g^i_{\gamma\alpha})d_{F^.}(\mu(e_\beta G^{i+1}_\alpha))$$
$$= -\sum\limits_{\gamma,\alpha}G^{i*}_\gamma
v(g^i_{\gamma\alpha})(f_\beta G^{i+1}_\alpha+\sum\limits_\delta
g^{i+1}_{\alpha\delta}\mu(e_\beta
G^{i+2}_{\delta}))=:II_1+II_2$$
plus a third one from
$\Hom(F^{i+1}, F^{i+1})$ (via pre-composing with $d_{F^.}$),
which equals
$$III:-\sum\limits_{\gamma,\alpha} G^{i+1*}_\alpha
v(g^{i+1}_{\alpha\delta})\mu(e_\beta G^{i+2}_\delta)\circ
d^i_{F^.}$$
$$=-\sum\limits_{\gamma,\alpha,\delta}G^{i*}_{\gamma}g^{i}_{\gamma\alpha}
v(g^{i+1}_{\alpha\delta})\mu(e_\beta G^{i+2}_\delta)$$
Now clearly $I+II_1=0$; and because of the derivation property
of $v$ and the characteristic syzygy property \refp{syzy}.,
we have $II_2+III=0$ as well.
Therefore the total is zero.
\par Now in the general case of $k\leq 0$, we consider
a section $$s=G^k_\beta\otimes v\in F^k\otimes T_P,$$
and we need
similarly to show the components of $d^2(s)$ in $\Hom(F^i,
F^{i+k+1}_+)$ vanish for all $i\leq 0$. Again there are 2
cases.
For $i=0$ there are the 2 contributions $$\sum\limits_\gamma
g^k_{\ba}v(f_\gamma) e^*_\gamma G^{k+1}_\alpha $$
and $$-\sum\limits_\gamma
g^k_{\ba}v(f_\gamma) e^*_\gamma G^{k+1}_\alpha $$
which cancel out. Then for $i<0$ there are again
3 contributions.
First via $F^{k+1}\otimes T_P$: $$I:
\sum\limits_{\gamma,\alpha,\eta}
g^k_{\beta\eta} v(g^i_{\gamma\alpha})G^{i*}_\gamma
\mu(G^{i+1}_\alpha G^{k+1}_\eta), $$ second via $\Hom(F^i,
F^{i+k})$:
$$II:-\sum\limits_{\gamma,\alpha\eta}G^{i*}_\gamma
v(g^i_{\gamma\alpha})d_{F^.}(\mu(G^k_\beta g^{i+1}_\alpha))$$
$$= -\sum\limits_{\gamma,\alpha\eta}G^{i*}_\gamma
v(g^i_{\gamma\alpha})(g^k_{\beta\eta}\mu(G^{k+1}_\eta
G^{i+1}_\alpha)+
\sum\limits_\delta
g^{i+1}_{\alpha\delta}\mu(G^k_\beta
G^{i+2}_{\delta}))=:II_1+II_2$$
and third via $\Hom(F^{i+1}, F^{i+k+1}_+)$:
$$III:-\sum\limits_{\gamma,\alpha} G^{i+1*}_\alpha
v(g^{i+1}_{\alpha\delta})\mu(G^k_\beta G^{i+2}_\delta)\circ
d^i_{F^.}$$
$$=-\sum\limits_{\gamma,\alpha,\delta}
G^{i*}_{\gamma}g^{i}_{\gamma\alpha}
v(g^{i+1}_{\alpha\delta})\mu(G^k_\beta G^{i+2}_\delta).$$
Again the 3 cancel out. This completes the proof that $\t_X$ is
a complex. Again, one can give a more abstract proof that
$d^2(v\otimes G)=0$, based on \refp{dTX-def-bracket}..
\par Note the fact that $\t_X$ is a complex, modulo knowing the
same for $\t_{X/P}$, could equivalently be formulated as
follows. Let $$J:F^.\otimes T_P\to T_{X/P}$$ be defined
by\begin{eqnarray}\label{Jdef}J^0(e_\alpha\otimes v)=\hskip
2in\\ \nonumber (f_\alpha v, \sum\limits_{i,\beta\gamma}
v(g^j_{\beta\gamma})G^{j*}_\beta\mu(e_\alpha
G^{j+1}_\gamma))\in T_P\oplus\n^0=\t_{X/P}^0, i=0\\
J^i(G^i\otimes v)=\hskip 2in\\ \nonumber
\sum\limits_{i,\beta\gamma}
v(g^j_{\beta\gamma})G^{j*}_\beta\mu(G^i G^{j+1}_\gamma))\in
\n^i=\t_{X/P}^i, i<0.
\end{eqnarray} (cf. \refp{J}.).
Then $J$ is a map of complexes, whose mapping cone is $\t_X$.
This what we have shown.

\subsubsection{ Actions}\ \par Next, we will define an action of
$\t_{X/P}$ on $ I.$ This action amounts to a pairing in the
derived
category
$$\t_{X/P}\times F^.\to F^.,$$ that extends to a morphism
\begin{equation}
\label{I action}\t_{X/P}\otimes F^.\to F^.\end{equation}
 (compatible with the differentials) . The pairing is defined
as
 follows: the component
  $$\Hom(F^i, F^j_+)\otimes F^i\to F^j$$ is the obvious
 map if $j\leq 0$ and 0 if $j=1$ (this is 'explained' by the
remark that via the inclusion $\n^.\to\gl(F^._+)$, an element
of $\n^.$ annihilates, hence preserves, the subcomplex
$F^1_+\subset F^._+$, hence induces an endomorphism of
$F^.=F^._+/F^1_+$); the component $$T_P\otimes F^i\to F^i$$
is defined by the postulation that it annihilate the standard
basis $(G^i_\beta)$ of $F^i$, thus
 $$v\otimes \sum\limits_\beta
 a_\beta G^i_\beta\mapsto \sum\limits_\beta
 v(a_\beta)G^i_\beta.$$ Now verifying that this yields a
pairing
 of complexes amounts to commutativity, for all $j<0$, of the
diagrams
 $$\begin{matrix}\Hom(F^i, F^j)\otimes F^i&\to& F^j\\
 \downarrow&&\downarrow\\ \Hom(F^i,
 F^{j+1})\otimes F^i&\to&F^{j+1}\end{matrix}$$
 (which is obvious); and, for all $i<0$, of
\begin{equation}\label{TF action}
 \begin{matrix}T_P\otimes
 F^i&\to&F^i\\ \downarrow&&\downarrow\\
 (\bigoplus\limits_j\Hom(F^j, F^{j+1})\otimes F^i)\oplus
 (T_P\otimes F^{i+1})&\to&
 F^{i+1}\end{matrix}\end{equation} To check commutativity of
 \refp{TF action}.,
consider a section $$s=v\otimes G^i_\beta \in T_P\otimes F^i.$$
 Then $s$ clearly goes to zero in $F^i$,
 hence in $F^{i+1}$ going clockwise. In
$ \bigoplus\limits_j\Hom(F^j, F^{j+1})\otimes F^i$, the
image of $s$  is $$-\sum\limits_{j,\gamma,\alpha} G^{i*}_\gamma
v(g^j_{\alpha\gamma})G^{j+1}_\gamma\otimes G^i_\beta$$ which
maps to
$$-\sum\limits_{\alpha}v(g^i_{\alpha\beta})G^{i+1}_\alpha\in
F^{i+1}$$
In $T_P\otimes F^{i+1}$, the image of $s$ is
$v\otimes\sum\limits_\alpha g^i_{\alpha\beta}G^{i+1}_\alpha$,
which maps to
$$\sum\limits_{\alpha}v(g^i_{\alpha\beta})G^{i+1}_\alpha\in
F^{i+1}.$$ The latter two add up to zero in $F^{i+1}$. Thus the
image of $s$ going counterclockwise is zero as well. This shows
the
diagram commutes, hence we have a pairing of complexes as
claimed.\par Note that a similar recipe also defines an action
of
$\t_{X/P}$ on $F^._+\sim A_X$, as well as on $A_P$ (the latter
via
the quotient $T_P$ of $\t_{X/P}$. Also, these actions are
derivations with respect to the $A_P$-module structure on $F^.$
and
$F^._+$: \beq\label{der-AP}\langle u, ab\rangle = a\langle
u,b\rangle + \langle u,a\rangle b, a\in A_P, b\in F^.,
F^._+.\eeq
Indeed because the action of any $u\in\n^.$ is $A_P$-linear and
such
$u$ annihilates $A_P$, it suffices to check this for $u\in T_P,
b=G^i_\alpha$, in which case it is immediate from the
definition.\par Note that it is not in general true that
\beq\label{der-AX}\act{u}{\mu(ab)}\equiv\mu(\act{u}{a}b)
\pm
\mu(a\act{u}{b})\mod dF^.\eeq where $dF^.$ is the image of the
map
$F^.[-1]\to F^._+$ induced by the differential $d$. For
example,
\refp{der-AX}. may fail if $u\in\n^2$. Nonetheless it is easy
to
see
that \refp{der-AX}. holds if $u$ is of degree $\leq 1$. In
particular, the induced pairing
$$H^0(\t_{X/P})\times A_P\to A_X$$ is a derivation on $A_P$,
hence
the pairing $$H^0(\t_{X/P})\times A_X\to A_X$$ is a derivation
on
$A_X$. Note that because $F^._+$ has cohomology only in degree
1
and \refp{der-AP}. holds for $a,b\in F^1_+$, the diagram
\begin{equation}\begin{matrix}&\t_{X/P}\times F^._+\times
F^._+&
\stackrel{\id\times\mu}{\to}&
\t_{X/P}\times F^._+&\\ &\downarrow&&\downarrow&\act{}{}\\
 &\t_{X/P}\times
F^._+&\stackrel{\act{}{}}{\to}&F^._+&\end{matrix}\end{equation}
commutes on cohomology.\par
 Next, we define in a similar fashion
an action
\begin{equation}\label{Lie ideal}\t_{X/P}\times ( I\otimes
T_P)\to
 I\otimes T_P\end{equation} via a pairing of complexes
\begin{equation}\label{IT action} \t_{X/P}\times (F^.\otimes
T_P)\to
F^.\otimes T_P\end{equation} Again, the key point is to define
$$T_P\times (F^i\otimes T_P)\to F^i\otimes T_P$$ by
$$v\times\sum aG^i_\alpha\otimes w_\alpha  \mapsto
\sum\limits_\alpha( aG^i_\alpha\otimes [v,w_\alpha]+
v(a)G^i_\alpha\otimes w_\alpha ), a\in A_P$$ (note that this is
compatible with our earlier definition setting
$v(G^i_\alpha)=0$). The pairing
$$\Hom(F^i, F^j)\times F^i\otimes T_P\to F^j\otimes T_P$$
is the obvious one (acting on the $F^i$ factor only). The
verification that this defines a pairing of complexes is again
essentially obvious on the $\n^.$ subcomplex of $\t_{X/P}$, and
it remains to check commutativity, for all $i<0$, of
\begin{equation}\label{TF action-IT}
 \begin{matrix}T_P\times
 F^i\otimes T_P&\to&F^i\otimes T_P\\ \downarrow&&\downarrow\\
 (\bigoplus\limits_j\Hom(F^j, F^{j+1})\otimes F^i\otimes
T_P)\oplus
 (T_P\times F^{i+1}\otimes T_P)&\to&
 F^{i+1}\otimes T_P\end{matrix}\end{equation} Indeed going
clockwise, an element \begin{equation}\label{clock}v\times
G^i_\alpha\otimes w\mapsto\sum\limits_\beta g^i_\ab
G^{i+1}_\beta [v,w].\end{equation} As for the counterclockwise
direction , first going downwards, $$v\times
G^i_\alpha\otimes w\mapsto(-\sum\limits_{\gamma\beta}
v(g^i_{\gamma\beta})G^{i*}_\gamma\otimes G^{i+1}_\beta \times
G^i_\alpha\otimes w, v\times\sum\limits_\beta
g^i_{\ab}G^{i+1}_\beta\otimes w)$$ Then going rightwards, the
first component maps to $$-\sum\limits_{\beta}
v(g^i_{\alpha\beta})G^{i+1}_\beta \otimes w $$
 and the second to
 $$\sum\limits_\beta g^i_\ab G^{i+1}_\beta
[v,w]+\sum\limits_{\beta} v(g^i_{\alpha\beta})G^{i+1}_\beta
\otimes w,$$ and the two add up to \refp{clock}..
\subsubsection{ $\t_{X/P}$ is a dgla} \ \par Next, we show that
$\t_{X/P}$ itself admits a bracket, making it a dgla. This
bracket is defined in the obvious way on element
pairs within $\n_{X/P}$ or
$T_P$. As for the cross terms $T_P\times\Hom(F^j, F^{i+j}_+)$,
 note first that
$\Hom(F^j, F^{i+j}))$ has a standard basis of the form
\begin{equation} h^{-j, i+j}_{\ab}=(G^j_\alpha)^*\otimes
G^{i+j}_\beta\end{equation} (i.e. the rank-1 homomorphism
taking
the
basis element $G^j_\alpha$ to $G^{i+j}_{\beta}$); with the
convention that when $i+j=1, (G^1_\beta)$ is the unique element
$1\in F^1_+$. We define
$$[v,h^{ij}_{\ab}]=0, \forall v\in T_P, \forall
i,j,\alpha,\beta$$
and extend by the derivation rule, i.e.
\begin{equation}\label{bracket_on_TX}
[v,\sum\limits_{\ab}c_{\ab}h^{ij}_{\ab}]=
\sum\limits_{\ab}v(c_{\ab})h^{ij}_{\ab}.\end{equation}\par Now,
the Jacobi
identity for this bracket amounts to 2
identities:
\begin{equation}\label{jac1}[[v_1,v_2],h]=[v_1,[v_2,h]]-[v_2,[v_1,h]],
\end{equation}\begin{equation}\label{jac2}
[v,[h_1,h_2]]=[[v,h_1],h_2]-[[v,h_2],h_1].\end{equation} Indeed
(\ref{jac1}) is obvious from the definition (both sides vanish
on
standard basis elements $h_{\ab}^{ij}$), while (\ref{jac2}) is
obvious from the bilinearity of the bracket on $h_1, h_2$ and
the
derivation property of the action of $v$ on scalars (which
again
just amounts to the fact that both sides vanish on the standard
basis). Thus we have a Lie bracket on $\t_{X/P}$.\par To
complete
the proof that $\t_{X/P}$ is a dgla,
it remains to verify the compatibility of the bracket with the
differential $d$.
Since this compatibility is already known within $\n$, we are
reduced to showing
\begin{equation}\label{tpx-dgla1}d[v,h]{=}[dv,h]+[v,dh], v\in
T_P,
h=h^{ij}_{\ab},\end{equation}\begin{equation}\label{tpx-dgla2}
d[v_1,v_2]=[dv_1, v_2]+[v_1,dv_2], v_1, v_2\in
T_P\end{equation}
where the LHS of \refp{tpx-dgla1}. is of course 0 by
definition.
Now we calculate, for $h=h^{-j, i+j}_{\ab}$: $$[dv, h] =
-\sum\limits_\gamma v(g^{i+j}_{\beta\gamma})
G^{j*}_\alpha\otimes G^{i+j+1}_\gamma
-(-1)^i\sum\limits_\delta
v(g^{j-1}_{\delta\alpha})G^{j-1*}_\delta
G_\beta^{i+j},$$ while
$$[v,dh]=[v, \sum\limits_\gamma g^{i+j}_{\beta\gamma}
G^{j*}_\alpha\otimes G^{i+j+1}_\gamma+
(-1)^i\sum\limits_\delta g^{j-1}_{\delta\alpha}G^{j-1*}_\delta
G_\beta^{i+j}]$$
$$=\sum\limits_\gamma v(g^{i+j}_{\beta\gamma})
G^{j*}_\alpha\otimes G^{i+j+1}_\gamma
+(-1)^i\sum\limits_\delta
v(g^{j-1}_{\delta\alpha})G^{j-1*}_\delta
G_\beta^{i+j},$$ which shows \refp{tpx-dgla1}.. Then
\refp{tpx-dgla2}. is equally easy: $$[v_1, dv_2]=[v_1,
-\sum\limits_{\ab} v_2(g^i_\ab)G^{i*}_\alpha G^{i+1}_\beta]=
-\sum\limits_{\ab} v_1(v_2(g^i_\ab))G^{i*}_\alpha
G^{i+1}_\beta,
$$$$[dv_1, v_2]=- [v_2, dv_1]= -[v_2,
-\sum\limits_{\ab} v_1(g^i_\ab)G^{i*}_\alpha
G^{i+1}_\beta]$$$$=
\sum\limits_{\ab} v_2(v_1(g^i_\ab))G^{i*}_\alpha G^{i+1}_\beta,
$$ and finally $$d[v_1, v_2]=
-\sum\limits_{\ab} [v_1,v_2](g^i_\ab)G^{i*}_\alpha
G^{i+1}_\beta,
$$ which yields \refp{tpx-dgla2}..
Thus we have shown that $\t_{P/X}$ is a dgla. Note that by
construction, $\n$ is a dg Lie ideal in $\t_{X/P}$ and we have
a
dgla exact sequence $$0\to\n\to\t_{P/X}\to T_P\to 0.$$\par
Next, we must check that the bracket just defined is compatible
with the actions of $\t_{X/P}$ on $F^.$ and on $F^.\otimes T_P$
defined in the pervious paragraph, i.e. that
\begin{eqnarray}\label{action-bracket} \langle [u_1,
u_2],a\rangle=\langle u_1,\langle u_2,a\rangle\rangle
-(-1)^{\deg(u_1)\deg(u_2)}\langle u_2,\langle
u_1,a\rangle\rangle,\\
\nonumber\forall u_1, u_2\in\t_{X/P}, a\in F^., F^.\otimes
T_P.\end{eqnarray}
Now in the case of the action on $F^.$, it is easy to see that
we may assume first that $a=G^i_\alpha$, and second that $u_1,
u_2$ are of the form $v$ or $h$ as above. Then in case $u_1$ or
$u_2=v$, both sides of \refp{action-bracket}. yield zero, while
if both $u_1, u_2\in\n\subset \gl(F^._+)$,
\refp{action-bracket}. holds by definition of the bracket on
$\gl$. The case $a=G^i_\alpha$ is similar.\par
Finally, we claim that the action of $\t_{X/P}$ in $F^.$is
compatible with the differential $d_{F^.}$, in other words,
that $d_{\t_{X/P}}$ maps under the action to commutator with
$d_{F^.}$, i.e
$$\act{d_{\t_{X/P}}(u)}{a}=d\act{u}{a}-(-1)^{\deg u}\act{u}{
d_{\F^.}(a)}, u\in\t_{X/P}, a\in F^..$$ But this is immediate
when $u\in \n$, and when $u\in T_P$ it follows directly from
\refp{dTXP-def-bracket}..
\subsubsection{ $\t_X$ is a dgla}\ \par
 Since we know $\t_X$ is a complex and brackets between
elements
of $\t_{X/P}$ and $F^.\otimes T_P$ have been defined previously
(cf. \refp{IT action}.), defining a dgla structure on $\t_X$
now amounts to
 defining the bracket on terms of the form $F^i\otimes T_P$,
and
verifying compatibility and Jacobi identity. To this end,
define
\begin{eqnarray}\label{ITP-dgla}
[a_\alpha G^i_\alpha v_1, a_\beta G^j_\beta
v_2]=\hskip 6cm \\ \nonumber
\sum\limits_\gamma a_\alpha v_1(a_\beta) \mu(G^i_\alpha
G^{j}_\beta) v_2-
\sum\limits_\delta a_\beta v_2(a_\alpha) \mu(G^{i}_\alpha
G^{j}_\beta) v_1 \\ \nonumber+a_\alpha a_\beta\mu(G^i_\alpha
G^j_\beta)[v_1, v_2]\\ \nonumber = \mu(G^i_\alpha G^j_\beta)
[a_\alpha v_1, a_\beta v_2]\end{eqnarray} To show $\t_X$ is a
dgla we first check compatibility of the bracket with the
differential , i.e. \begin{equation}\label{d[]} d[a_\alpha
G^i_\alpha v_1, a_\beta G^j_\beta
v_2]=[d(a_\alpha G^i_\alpha v_1), a_\beta G^j_\beta
v_2]+(-1)^i[a_\alpha G^i_\alpha v_1, d(a_\beta G^j_\beta
v_2)].\end{equation}
It is easy to see that we may assume $a_\alpha=a_\beta=1.$
In this case the LHS of \refp{d[]}. is just
$$d\mu(G^i_\alpha G^j_\beta)[v_1,v_2]+\mu(\mu(G^i_\alpha
G^j_\beta)d[v_1, v_2])=$$
$$\mu(dG^i_\alpha G^j_\beta)[v_1,v_2]-(-1)^i\mu(G^i_\alpha
dG^j_\beta)[v_1,v_2]+\mu(\mu(G^i_\alpha
G^j_\beta)d[v_1,v_2])$$ On the other hand the first bracket of
the RHS yields $$[d(G^i_\alpha) v_1+\mu(G^i_\alpha dv_1),
G^j_\beta v_2]=$$
$$ \sum\limits_\gamma v_2(g^i_{\alpha\gamma})\mu(
G^{i+1}_\gamma G^j_\beta) v_1 + \mu(dG^i_\alpha
G^j_\beta)[v_1,v_2]+$$
$$\sum\limits_\delta v_1(g^j_{\beta\delta})\mu(G^i_\alpha
G^{j+1}_\delta)v_2+\mu(\mu(G^i_\alpha G^j_\beta)[dv_1,v_2])
$$

while the second bracket yields
$$[G^i_\alpha v_1, d(G^j_\beta)v_2+\mu(G^j_\beta dv_2)]=$$
$$(-1)^{i+1}\sum\limits_\delta v_1(g^j_{\beta\delta})
\mu(G^i_\alpha G^{j+1}_\delta)v_2+\mu(G^i_\alpha
dG^j_\beta)[v_1,v_2]+$$
$$(-1)^{i+1}\sum\limits_\gamma v_2(g^i_{\alpha\gamma})
\mu(G^j_\beta G^{i+1}_\gamma) v_1
+(-1)^i\mu(\mu(G^i_\alpha G^j_\beta)[v_1,dv_2]))  $$
Therefore \refp{d[]}. holds. \par It remains to verify the
Jacobi
identity, in the form\begin{eqnarray}\label{jac-IT}
[[G^i_\alpha
v_1, G^j_\beta v_2], G^k_\gamma v_3]= \\ \nonumber[G^i_\alpha
v_1, [G^j_\beta v_2, G^k_\gamma
v_3]]-(-1)^{(i+1)(j+1)}[G^j_\beta v_2, [G^i_\alpha v_1,
G^k_\gamma v_3]]\end{eqnarray} As remarked earlier, we may
assume that $\mu(G^i_\alpha G^j_\beta)$ is a standard basis
element of $F^{i+j-1}$, of the form $G^{i+j-1}_\delta$ for some
$\delta$, and similarly for other $\mu$ products. In this
case, \refp{jac-IT}. follows directly from the definition
\refp{ITP-dgla}.. This finally completes the proof the $\t_X$
is
a dgla.
\par As for the extension of the action of $\t_{X/P}$ on
$F^._+$
to an action of $\t_X$, this is defined by setting
$$\act{G^i_\alpha v}{ aG^j_\beta}= v(a)\mu(G^i_\alpha
G^j_\beta), i\leq 0, j\leq 1$$ (note that $G^i_\alpha\otimes v$
has degree $i-1$ in $\t_X$, so both sides have the same degree,
viz. $i+j-1$, as required).
The same formula also defines an action
\beq\label{tx-f+-action}\t_X\times F^._+\to F^._+.\eeq
 It is easy to check that this is a pairing of complexes, and
that it is compatible with the bracket on $\t_X$, and that all
elements of degree $\leq 1$ in $\t_X$ acts on $F^._+$ as
derivations $\mod$ the image of $ F^.$.  Therefore, we get a
derivation  action
\begin{equation} H^0(\t_X)\times A_X\to A_X\end{equation} which
by its very construction is compatible with that of $T_P$ on
$ A_P$ and $\t_{X/P}$ on $ I$.\par We claim next that $\t_X$ is
an
$ A_X$-module, in the sense that there is (in the derived
category) a multiplication pairing
\beq\label{module}  A_X\times\t_X\to\t_X,\eeq, or to be
precise,
a map of complexes with appropriate 'action' properties \beq
F^._+\times \t_X\to \t_X[-1].\eeq  This is simply defined by
postmultiplication:
$$\act{G^i_\alpha}{G^j_\beta\otimes v}= \mu(G^i_\alpha
G^j_\beta)\otimes v,$$
$$\act{G^i_\alpha}{ G^j_\beta\otimes G^{k*}_\gamma}=
\mu(G^i_\alpha G^j_\beta)\otimes G^{k*}_\gamma.$$
 Note that the $ A_X$-module
structure (\ref{module}) shows that the cohomology sheaves of
$\t_X$
are $ A_X$-modules compatibly with their $ A_P$-module
structure,
hence clearly $ A_X$-coherent. In particular $\t_X$ is
equivalent to
its sheaf-theoretic restriction on $X$. Since $\t_X$ is a
bounded
complex with $ A_X$-coherent cohomology, it defines an element
of
the bounded, coherent derived category $D^b_c(X).$\par We note
that
even though in the above construction we used the 'full dgla'
model
for the normal atom $\n_{X/P}$, leading to a dgla model for
$\t_X$,
we could instead have used the reduced model $\n^{\red}$ (see
\S\ref{rednormalsec})
for $\n$, thus leading to a Lie atom $\t_X^{\red}$ weakly
equivalent
to $\t_X$.
\subsubsection{Dependence on choices}\label{dep-on-choice}\ \nl
The above construction of the tangent atom $\t_X$ of an affine
scheme $X$ depended, in addition to a choice of affine
embedding
$X\to P$, also on a choice of free resolution of $I=I_{X/P}$.
We claim next that, still fixing the embedding   $X\to P$, the
Lie atom $\t_X$ is independent of the resolution \emph{up to
Lie
quasi-isomorphism}, i.e. a composition of Lie atom
homomorphisms
inducing an isomorphism on cohomology and their inverses. This
is just a trivial variation of the usual statement on
independence of resolution. The precise statement is as
follows:
\begin{lem}\label{lem123}. Let $^1F^., ^2F^.$ be free
resolutions of the ideal
of the closed subscheme $X\to P$ and $^1\t_X, ^2\t_X$ the
associated Lie atom. Then there exists a third such resolution
$^3F^.$  with associated atom $^3\t_X$,  together with direct
Lie quasi-isomorphisms to the associated Lie atom
\beq\label{tx123}^1\t_X\to^3\t_X\leftarrow
^2\t_X\eeq\end{lem}\begin{proof}
As is well known, there exists a resolution $^3F^.$ of $I$
together with maps $$^ip:^3F^.\to ^iF^., ^iq:^iF^.\to ^3F^.$$
such that $^ip^iq$ is the identity and $^iq^ip$ is homotopic to
the identity for $i=1,2$. Moreover we may assume $^iq$ takes a
standard generator $G^j_\alpha$ to a standard generator. Then
we
get maps $$^i\n=\Hom(^iF^., ^iF^._+)\to ^3\n, i=1,2,$$
$$h\mapsto ^iq\circ h\circ ^ip.$$ These are clearly Lie
quasi-isomorphisms. By construction,
 these extend to maps as in \refp{tx123}..
\end{proof} Next we take up the question of dependence of
$\t_X$
on the affine embedding $X\to P$. To make this dependence
explicit, we will henceforth denote by $\t_X(P)$ an explicit
representative of the Lie atom denoted formerly by $\t_X$,
based
on a choice of resolution $F^.$, and uniquely determined up to
Lie quasi-isomorphism.
 What we claim is that the weak
equivalence class of $\t_X(P)$ is independent of the embedding;
more
precisely, given two embeddings $X\to P, X\to Q$, there is a
Lie
atom $\t$ and (direct) weak equivalences
\beq\label{emb-compare1}\t_X(P)\to\t, \t_X(Q)\to\t\eeq (we
recall that a
weak equivalence is not necessarily invertible even in the
derived
category; a direct weak equivalence is one given by a morphism
of
complexes). In fact, we claim that can take
\beq\label{emb-compare2}\t=\t_X(P\times Q)\eeq where $X\to
P\times
Q$ embeds via the diagonal $\Delta_X\subset X\times X$.  What
we need to do then is construct each of the weak equivalences
(\ref{emb-compare1}).\par To this end set $R=P\times Q$. Note
that the embedding $X\to R$ lifts to $P\to R$, which is a cross
section of projection $R\to P$. We first contrcut a Lie
homomorphism \beq\label{txp to txr}\t_X(P)\to \t_X(R).\eeq This
is done as follows.
Let $F^._P$ be a free resolution of
$ I_{X/P}$ 
Then we can construct a free resolution $F^._R$
of $ I_{X/R}$ such that each $F^i_R$ contains $F^i_P\otimes
A_R$
as
a direct summand. In fact, since  the embedding $P\to R$ is
automatically regular, we may assume that $$F^i_R=F^i_P\otimes
A_R\oplus F^i_Q$$ where $F^._Q$ is a Koszul resolution of
$I_{P/R}$. Set $$S^.=\txp\otimes A_R$$ which has an obvious
structure of Lie atom. Then as in the proof of  Lemma
\ref{lem123} above,
we get a Lie homomorphism $$\n_{X/P}\otimes
A_R=\Hom^._R(F^._P\otimes A_R, F^._{P+}\otimes A_R)\to
\n_{X/R}=\Hom^.(F^._R, F^._{R+})$$
which extends to a Lie homomorphism $S\to\t_X(R)$.
Combined with the obvious homomorphism $\txp\to S$, this gives
the
desired homomorphism \refp{txp to txr}.. Note that this
homomorphism can be represented schematically as follows
\begin{equation}\label{tp-s-r}\begin{matrix}T_P\otimes
F^._{P+}&\to&\n_{X/P}\\
\downarrow && \downarrow \\
T_P\otimes F^._{P+}\otimes A_R & \to & \n_{X/P}\otimes A_R\\
\downarrow&& \downarrow\\
T_R\otimes F^._{R+}&\to&\n_{X/R}\end{matrix}\end{equation}
Note that since $F^._P\otimes A_R$ is a resolution of
$I_{X\times Q\to P\times Q}$, we can represent the above
diagram
\refp{tp-s-r}. as follows
\beq\label{tpsr}\begin{matrix}T_P\otimes
F^._{P+}&\to&\Hom^.(F^._P, A_X)\\
\downarrow&& \downarrow \\
T_P\otimes F^._{P+}\otimes A_R & \to & \Hom^.(F^._P, A_{X\times
Q})\\
&& \downarrow\\
\downarrow&&\Hom^.(F^._P, A_{X\times Q})\otimes A_X\\
&&\downarrow\\
T_R\otimes F^._{R+}&\to&\Hom^.(F^._R,
A_X)\end{matrix}\end{equation} Now the composite of the left
arrows is clearly an injection with cokernel $T_Q\otimes A_X$,
while the composite of the 2 upper right arrows is an
isomorphism. On the other hand as in \S\ref{rednormalsec} we
have that $\Hom(F^._Q, A_X)=\n_{P/R}\otimes A_X$ is weakly
equivalent to the restricted reduced normal atom
$$\n^{\red}_{Q/R}\otimes A_X=N_{P/R}[-1]\otimes A_X\simeq
T_Q[-1]\otimes A_X.$$ Similarly, we have an exact sequence
\begin{figure}
\end{figure}
$$\n_{X/P}^{\red}\to\n_{X/R}^{\red}\to T_Q[-1]\otimes A_X\to$$
 which means that, up to weak equivalence, the cokernel of the
bottom right arrow is $T_Q\otimes A_X$. Therefore the map from
the top row of \refp{tpsr}. is a weak equivalence, as claimed.

\subsubsection{Functoriality}\label{functoriality}\
 \par We claim next, still in
the embedded case, that a morphism
\beq\label{fXY}f:X\to Y\eeq of
affine schemes induces a map
\beq\label{df} df:\t_X\to f^*\t_Y.\eeq  A precise formulation
is as follows. First, note that given
affine embeddings $$X\to P, Y\to Q,$$ a morphism $f$ as
in \refp{fXY}. extends to a morphism
$f:P\to Q$. Moreover since $\t_Y(Q)$ is a complex of finite free
$A_Q$-modules, $f^*\t_Y=f^!\t_Y$ refers unambiguously to the
termwise pullback, which is, in fact, a Lie atom (cf.
\cite{atom}).\par Then the precise assertion is that  there exists another
weak-equivalence representatives $\t_X(P_1)$
of
$\t_X(P)$,  
 an extension of $f$ to a morphism
$\tilde{f}:P_1\to Q$
with a map of Lie atoms \beq\label{df1}df:\t_X(P_1)\to
\tilde{f}^*\t_Y(Q)\eeq
(here we mean an 'actual' map of complexes, without inverting any weak
equivalences). Moreover $\tilde{f}^*\t_Y(Q)$ admits a structure of
a $\t_X$-module such that $df$ is $\t_X$-linear. Note that this
would make the
mapping cone of (\ref{df1})
itself into a Lie atom, which (compare \cite{atom}) may be
denoted
$\n_{f/Y}[-1]$). \par Now the construction of \refp{df1}. is the
following.
Identifying $X$ with the graph of $f$, we have a diagram
\beq\begin{matrix} X&\to& P\times
Y&\to&Y\\&\searrow&\downarrow&&\downarrow\\ &&P\times Q&\to&
Q\end{matrix}\eeq which induces a diagram
\beq\label{df2}\begin{matrix}
T_{P\times Q}\otimes A_X&\to&\n_{X/P\times Q}\\
\downarrow&&\downarrow\\
T_Q\otimes A_X&\to&f^*\n_{Y/Q}\end{matrix}\eeq where the right
vertical arrow comes from the fact that a resolution of $Y\to
Q$
pulls back to a resolution of $P\times Y\to P\times Q$, hence
can be
assumed so be a subcomplex, and termwise a direct summand, of a
resolution of $X\to P\times Q$. Then (\ref{df2}) yields a map
as
in
(\ref{df1}) with $P_1=P\times Q$, and it is straightforward to check that this map is
indeed $A_X$-linear.\par
Finally, we claim that the 'functoriality map' $df$ is indeed
functorial, in the sense that given morphisms of affine schemes
$$X\stackrel{f}{\to}Y\stackrel{g}{\to} Z$$
then commutes\beq \begin{matrix}
\t_X&&\ \ \ \ \stackrel{df}{\longrightarrow}\ \ \ &f^*\t_Y\\
&&\\
&\searrow&d(g\circ f)&\downarrow f^*dg\\
&&\\
&&&f^*g^*\t_Z.\end{matrix}\eeq
This can be made precise (via explicit representatives), and
proven very similarly to the above.
 This finally completes the proof of Theorem
\ref{affine tx}.\end{proof}

\begin{cor}\label{smoothaffine}
 If $X$ is a smooth affine variety, then
the weak equivalence class of $\t_X$ coincides with that of the
usual tangent module $T_X=\Der(A_X)$.\end{cor}\begin{proof} Any
affine embedding $X\to P$ is regular, so we may use for the
normal
atom $\n_{X/P}$ the shifted normal module $N_{X/P}[-1]$. In
this
case $\t_X$ may be identified with the kernel of the surjection
$$T_P\otimes A_X\to N_{X/P},$$ which, as is well known, is just
$T_X$.\end{proof}
\subsection{Sheafification} Let $\t_X(P)\st$ be the termwise
sheafification of $\t_X(P)$, as complex of sheaves of (free)
$\O_P$-modules on $P$. Clearly, $\t_X(P)\st$ is acyclic off $X$,
hence is quasi-isomorphic to its (topological) restriction on
$X$, i.e. the natural map
$$\t_X(P)\st\to i_{X*}i_X\inv\t_X(P)\st$$ is a quasi
isomorphism. Hence, we may consider $\t_X(P)\st$ as a complex of
sheaves on $X$ (though the sheaves are not $\O_X$-modules).
We have
\beq \t_{X/P}=\Gamma(P, \t_{X}(P)\st)\sim \Gamma(X,
i_X\inv\t_X(P)\st)\eeq

As a special case of functoriality, consider a function $0\neq
f\in A_P$ and the associated principal affine open
$X_f\stackrel{i}{\hookrightarrow} X$, which fits in a diagram
$$\begin{matrix} X_f&\stackrel{i}{\hookrightarrow}& X\\
\downarrow&&\downarrow\\
P\times\A^1&\stackrel{\pi_P}{\to}&P
\end{matrix}$$ Then it is easy to see that the functoriality
map
$$\t_{X_f}(P\times\A^1)\to i^!\t_X(P)$$
is a quasi-isomorphism, as is the natural map
\beq \t_{X_f}(P\times\A^1)\to \t_X(P)_{f}\eeq where the latter
refers to termwise localization.
\subsection{Maps of
affine schemes: tangent dgla}  Let
$$f:X\to Y$$ be a mapping of affine schemes. It is natural to
consider tangent data to $f$, that is, a compatible set of
tangent data to $X,Y$ and $f$. This can be done as follows.
\par Given affine embeddings
$X\to P, Y\to Q$, $f$ can be extended to a map $P\to Q$.
Replacing
$X\to P$ by the graph embedding $X\to P\times Q$, we may assume
$P\to Q$ is a product projection. Then we have an injection
$I_{Y,Q}\to I_{X,P}$ which extends to the free resolutions
$F^._Y\to
F^._X$, and we may moreover assume that each $F^i_Y\otimes
A_P\to
F^i_X$ is a direct summand inclusion. We can identify the
functor
$f^!$ on complexes with $f^!\cdot =\cdot\otimes_{A_Y}
F^._{X+}.$
Then the complex $f^!(\n_{Y/Q})$ can be represented by
$$\hom^._{A_Q}(F^._Y, F^._{+X})=\hom_{A_P}^.(F^._Y\otimes A_P,
F^._{+X})$$ and there are maps
\begin{equation}\label{nf}\n_{X/P}\to
f^!(\n_{Y/Q})\leftarrow \n_{Y/Q}\end{equation} The mapping cone
of
\refp{nf}. can be represented by the sub-dgla of
$\n_{X/P}\oplus\n_{Y/Q}$ consisting of pairs $(a^.,b^.)$ such
that
$a^.$ vanishes on the subcomplex $F^._Y\otimes A_P\subset
F^._X$. We
denote this mapping cone by $\n_f$ or, more properly,
$\n_{f,P,Q}$, and
refer to it as the \emph{normal dgla} of $f$.\par Next,
proceeding
as in the case of schemes, we can construct a suitable
representative of the mapping cone $K$ of
$$T_P\otimes F^._{+X}\to T_Q\otimes F^._{+X}\leftarrow
T_Q\otimes
F^._{+Y},$$ together with a map of $K$ to $\n_f$, so that the
mapping cone of $K\to\n_f$ is a dgla, called the \emph{tangent
dgla}
to $f$ and denoted $\t_f$ or more properly, $\t_f(P,Q).$ By
construction, $\t_f(P,Q)$ is the mapping cone of
\begin{equation}\label{tfseq}\t_X(P)\oplus\t_Y(Q)\to
f^!\t_Y(Q).\end{equation}
\subsection{Reduced tangent algebra}
As was the case for the normal algebra, the tangent algebra
$\t_X(P)$ admits a reduced, weakly equivalent version
$^\red\t_X(P)$ which is convenient for applications because it
has a smaller $H^2$ (in fact, the $H^2$ will vanish when $X$ is
a locally complete intersection, which
is not the case for $\t_X$ itself).
This dgla is closely related
to, essentially a dual version of, the
Grothendieck-Lichtenbaum-Schlessinger \emph{cotangent complex}
of $X$ (cf. \cite{lich}).\par
\begin{thm} Given an affine embedding $X\to P$, there is a dgla
$\tred_X(P)$, together with a direct weak equivalence
$$\tred_X(P)\to\t_X(P),$$ such that for all $i>1$,
$H^i(\tred_X(P))$ is supported on the locus of non-lci points of
$X$. Moreover, the dgla quasi-isomorphism class of $\tred_X(P)$
depends only on the isomorphism class of $X$ as scheme over
$\C$.\end{thm}
\begin{proof}
Recall that the reduced normal algebra $^\red\n=^\red\n_{X/P}$
admits an injective dgla map
$$\lambda:\nred\to\n'\stackrel{\id+\pi}{\longrightarrow}\n$$
where $\nred\to\n$ is a direct sum inclusion (see
\refp{n'split}.). Let
$$\rho:\n\stackrel{\id-\pi}{\to}\n'\to\nred$$
be the obvious left inverse to $\lambda$
(where $\n'\to\nred$ is the projection). Recall (see
\refp{Jdef}.) that $\t_{X}(P)$ could be defined as a mapping
cone of $$J:T_P\otimes F^._+\to\n.$$ Now by  the definition of
$J$, we have, for $v\otimes G\in T_P\otimes F^i_+$,
$$J(v\otimes G)(e_\alpha)=v(f_\alpha)G,$$
$$J(v\otimes G)(s_\ab)=v(f_\alpha)\mu(e_\beta\otimes G).$$
Therefore the component of $J(v\otimes G)$ in $\Hom(F\inv_0,
F^0)$ is equal to $\pi$ of its component in $\Hom(F^0, A_P)$,
i.e.
$$J(v\otimes G)|_{F\inv_0}=\pi(J(v\otimes G)|_{F^0}).$$
  Thus it follows that the image of $J$ is contained in that of
$\id+\pi$, hence $J$ factors through a map
$$^\red J=\rho\circ J:T_P\otimes F^._+\to \nred,$$
i.e. $J=\lambda\circ\ ^\red J.$
Set
\beq\label{tred-def}\tred_X(P)={\mathrm{mapping\ cone}}(^\red
J).\eeq
Then the fact that $\nred\to\n$ is an injective Lie homomorphism
and a weq and that $\t_X$ is a dgla implies that $\tred_X(P)$ is
a dgla, weakley equivalent to $\t_X$.
\par
Now to complete the proof, it remains
to prove the lat statement, i.e. independence of choices.
To this end, we modify the argument of subsection
\ref{dep-on-choice} as follows, using notation as there.
Choosing standard resolutions of $X\to P, X\to R$, we obtain
such
for $X\to R$, whence compatible complexes $F^._{1,P}, F^._{1,R}$
etc. and then we have  dgla homomorphisms
$$\nred_{X/P}\to \nred_{X/P}\otimes A_R\to \nred_{X/R}.$$
Then we have a commutative diagram
\beq\begin{matrix} T_P\otimes F^._{P+}&\to&\nred_{X/P}\\
\downarrow&&\downarrow\\
T_R\otimes F_{R+}&\to&\nred_{X/R}\end{matrix}\eeq where the
left column (as mapping cone) is clearly naturally
quasi-isomorphic to $T_Q\otimes A_X[-1]$. On the other hand the
discussion in subsection \ref{dep-on-choice}, together with
Corollary \ref{nred-xyz} show that the same is true of the right
column. Thus, the horizontal maps together give a dgla
quasi-isomorphism $$\tred_X(P)\to\tred_X(R).$$ This yields our
conclusion.
\end{proof}
Since $\t_X$ is defined to be the weak equivalence class of
$\t_X(P)$, we may take $\tred_X(P)$ as a representative of
$\t_X$; we call representatives obtained in this way
\emph{reduced}.
\subsection{Comparison}
The complex $\t_X$ is a bit complicated, therefore in applications it is important to be able to compare it to
something simpler, like a module or the Ext-dual of a module.
The following result, based on Ischebeck's Theorem (see
\cite{isch} or \cite{mat}, p. 104) is useful in this connection.\par
First a definition. For any algebraic scheme $X$, we denote by $\NCI_X\subset X$ the (closed) locus of point $x\in X$ where $X$ is not locally a complete intersection.
\begin{thm}\label{compare-to-ext} Let $X\to P$ be an affine scheme.Then \begin{enumerate}
\item there is a natural map \beq\label{tau}\tau_X:\Hom^.(\Omega_X, A_X)\to \tred_X(P);\eeq
\item the induced map on cohomology
\beq H^i(\tau_X):\ext^i(\Omega_X, A_X)\to H^i(\tred_X(P))\eeq is always bijective for $i=0$ and injective for $i=1$;
moreover, $H^1(\tau_X)$ is bijective and
$H^2(\tau_X)$ is injective provided $X$ is reduced (more generally, a generic locally complete intersection with no embedded points);
\item $\tau_X$ induces an isomorphism on cohomology in degrees $\leq k$ provided, for all $x\in X$, that
\beq\label{depth}{\mathrm{depth}}_xA_X-\dim_x\NCI_X\geq k.\eeq
\item In particular, if $X$ is Cohen-Macaulay and lci in codimension 1, then $\tau_X$ is a cohomological isomorphism in degrees $\leq 2,$ i.e.
\beq\label{extomega} H^i(\tau_X):\ext^i_X(\Omega_X, A_X)\simeq H^i(\tred_X(P)), i\leq 2.\eeq\end{enumerate}
\end{thm}
\begin{proof}
Consider the following complex in degrees $\leq 1:$
$$G^.:\ \ ...\to F\inv_1\otimes A_X\to F^0\otimes A_X\to\Omega_P\otimes A_X.$$
This is a complex of free $A_X$-modules, whose cohomology in degree 1 is clearly $\Omega_X$, whence a morphism
$$\check{\tau}_X:G^.\to \Omega_X.$$ Locally at any lci point $x\notin\NCI_X$, the resolution $F^.$ is equivalent to a Koszul resolution, and therefore the complex $G^.$ is exact at $x$
in degrees $<1$.
Thus for $i<1$, $H^i(G)$ is supported on $\NCI_X$.
Now the complex $\Hom^._X(G^., A_X)$ is quasi-isomorphic to
$\tred_X(P)[1] $, therefore by dualizing $\check{\tau}_X$ we get $\tau_X$. Moreover we have a standard spectral sequence
$$E_2^{p.q}=\ext^p_X(H^q(G^.), A_X)\Rightarrow H^{p+q+1}(\tred_X(P)).$$
From this spectral sequence the assertions of (ii) follow easily; for example the kernel of $H^2(\tau_X)$ comes from $\ext^0(H^{-1}(G), A_X)$, hence vanishes provided the support of $H^{-1}(G)$, i.e. $\NCI_X$, is a proper subscheme and $X$ has no embedded points.\par Finally for part (iii) we use
 Ischebeck's Theorem which says that the $\ext^p$ vanishes provided $q<1$ and $$p<{\mathrm{depth{A_X}}-\dim{\NCI_X}}.$$
Therefore $\tau_X$ is an isomorphism in the range as claimed.
\end{proof} Note that the condition of part (iii) of the theorem applies whenever $X$ is either normal and Cohen-Macaulay, or a locally complete intersection curve. For deformation theory, it is the cohomology in degrees $\leq 2$ of $\t_X$ that matters.
\subsection{Example: a limit curve} We consider an example
extending ones which occur in the
local study of limits of covering maps of curves (\cite{HM}, extending earlier work by Beauville) and families of nodal plane curves (Severi varieties)\cite{sing}. Let $C\subset \A^2_{u,v}$ be the curve with equation $uv$. It has components $C_1=\A^1_u, C_2=\A^1_v.$
 Its tangent algebra is
\beq \t_C: T_{\A^2}\otimes F^._{C,\A^2,+}\to\n_C=\Hom(F^._{C,\A^2}, F^._{C,\A^2,+}),\eeq
and may be written symbolically as
\beq A_C<\del_v,\del_u>\stackrel{uv}{\to} A_C,\eeq
or more simply as
\beq 2A_C\stackrel{(u,v)}{\longrightarrow} A_C.\eeq
Note that a basis  for $\n_C$ is given by $\delta_C=[uv]^*,$
i.e. the map that takes the value 1 on the resolution element $[uv]$ corresponding to the ideal generator $uv.$\par
Now let $X\subset\A^{n+2}_{x,y,z}=:Q, n\geq 0$ be the subscheme with equation $xy$, i.e. a copy of $C\times\A^n,$ with components
$X_1=\A^{n+1}_{x,\underline{z}}, X_2=\A^{n+1}_{y,\underline{z}}.$ The tangent algebra $\t_X$
is similarly given symbolically by the complex
\beq (n+2)A_X\stackrel{(x,y)}{\longrightarrow}A_X,\eeq
where the last summand is identified with $\n_{X/\A^{n+2}}$ via $\delta_X=[xy]^*.$
Then $H^1(\t_C)$ and $H^1(\t_X)$ are both 1-dimensional and there is no higher cohomology, which means that the universal deformation of $C$ and $X$ are respectively given by
\beq uv-t, xy-s.\eeq\par
Now consider the map
\beq f:C\to X,\ \  (u,v)\mapsto (u^m, v^m, f_1^u+f_1^v,...f_n^u+f_n^v)\eeq where the $f_i^u, f_i^v$ are unramified and without constant term. For example,\begin{itemize}\item if $n=0$, this is
an $m$-fold branched cover on each component;\item if $n=1$,
we may assume $f_1^u=u, f_1^v=v$ and then
this map sends $C_1, C_2$ to the curves in $X_1, X_2$ with respective equations $x-z^m, y-z^m$, while $f(C)$ in total has equation $x+y-z^m.$\end{itemize} Then the Lie atom $f^!\t_X$ measures deformations of $(f,X)$ with $C$ fixed, together with a trivialization of
the corresponding deformation of $X$. This atom can be written symbolically as
\beq\label{ft_X} (n+2)A_C\stackrel{(u^m, v^m)}{\longrightarrow} A_C\eeq
in which it is convenient to identify $A_C$ as $F^._{C,\A^{n+4}, +}[1]$, where $C\to \A^{n+4}_{u,v,x,y,\underline{z}}=:P$ is the graph embedding, with equations
\beq\label{eq-C-in-P}uv, xy, x-u^m, y-v^m,  z_i-f_i^u-f_i^v, i=1,...,n.\eeq  Thus we may rewrite \refp {ft_X}. symbolically as
\beq {n+2}F^._{C,P}\to {n+2}A_P\oplus F^._{C,P}\to A_P\eeq
where the last summand is generated by $\delta_X.$
It is easy to check that
\begin{eqnarray} H^0(f^!\t_X)=\{c_x(u)\del_x+c_y(v)\del_y+c_{\ul{z}}(u,v)\del_{\ul{z}}:\\
\nonumber c_x, c_y, c_{\ul{z}}\in A_C, c_x(0)=c_y(0)=0\},\\
H^1(f^!\t_X)=(A_C/(u^m, v^m)).\eps_X.\end{eqnarray}
Moreover, the bracket (i.e. the obstruction) vanishes on $H^0(f^!\t_f)$. Thus, $\Deff(f^!\t_X)$ is a smooth, albeit
infinite-dimensional space; in fact, it is quite elementary that $\Deff(f^!\t_X)$ can be identified with the space of pairs $(f_u, f_v)$ where $f_u, f_v$ are deformations, respectively,
of the restrictions of $f$ to the $u$ and $v$- axes, considered as maps to  $(x,\ul z)$ (resp. $(y,\ul z)$)- space, that map the origin to the same point on the $z$-hyperplane; and as such, this space is clearly smooth.
\par
Now consider the Lie atom $\t_f$, classifying general deformations of $f$. This can be written in the form
\begin{eqnarray}\label{t_f} A_C<\del_v, \del_v>\oplus
A_X<\del_x, \del_y, \del_{\ul{z}}>\oplus F^._{C,P}<\del_x, \del_y, \del_{\ul{z}}>
\\ \longrightarrow  A_C\delta_C\oplus A_X\delta_X\oplus A_P<\del_x, \del_y, \del_{\ul{z}}>\oplus F^._{C,P}\delta_X\\ \stackrel
{((uv)^{m-1},\id, (u^m,v^m,0),  \eps)}{\longrightarrow}A_P\delta_X
\end{eqnarray} where, wherever necessary, we replace $A_C$ by $F^._{C,P+}$ and $A_X$ by $F^._{X,Q+}$. Here each complex term, such as $F^._{C,P}$, should be viewed as extending leftward from the indicated position. Note that the last map reflects the fact that the differential of $\t_f$
maps $\delta_C=[uv]^*$ to \beq\label{df-delc}df(\delta_C)=(uv)^{m-1}[xy]^*=(uv)^{m-1}\delta_X
\eeq\par
As for brackets, note that
\beq\label{[delc, delx]}[\delta_C, [uv]\delta_X]=(uv)^{m-1}\delta_X\in A_P\delta_X=\t_f^2;\eeq otherwise, the obvious generators of $\t_f^1$ have zero bracket.
In particular we may, in computing the Jacobi-Bernoulli cohomology of $\t_f$, replace $\t^2_f$ by its quotient by the remaining generators of $F^0_C$
besides $uv$ (cf. \refp{eq-C-in-P}.), while eliminating the corresponding summands form $F^0_{C,P}\delta_X$.\par
Now consider a prospective special multiplicative cocycle
$\epsilon=\eps(\phi, \psi)$
 for
$\t_f$, and write $\phi, \psi$ in the form
$$\phi= (a\delta_C, b\delta_X)\in\t_C^1\oplus\t_X^1, \psi=(c_x\del_x+c_y\del_y, d[uv]\delta_X)\in \t_X^0\oplus F^0_{C,P}\otimes\t_X^1\subset (f^!\t_X)^0,$$
where $a,b$ are constants, $d=\sum\limits_{i,j<m}d_{i,j}u^iv^j[uv]$,
and by \refp{df-delc}. we may assume $d_{m-1,m-1}=0$.
Now the Jacobi-Bernoulli cocycle condition on $\eps$ reads
\beq\label{cocyc}i(\phi)+d(\psi)+\sum\limits_{i=1}^\infty B_i\ad(\psi)^i(\phi)=0 \eeq
where the $B_i$ are the Bernoulli numbers. Here all higher
iterated adjoints (with $i>1$) vanish for reasons of degree. Thus
\refp{cocyc}.
 yields
\beq\label{cocycle-df} a(uv)^{m-1}+b+c_xu^m+c_yv^m+duv+B_1ad=0\eeq
where, of course, $B_1$ is the 1st Bernoulli number $B_1=-1/2.$
Now, working modulo $(u^m, v^m)$ we see firstly that $d_{i,j}=0$ for $i\neq j$. Then, reading off coefficients of $(uv)^i$, we see that
\begin{eqnarray} b=ad_0/2\\
d_0=ad_1/2\\
...\ \  d_i=ad_{i+1}/2, i\leq m-3,\\ d_{m-2}=-a\end{eqnarray}
In particular, \beq\label{a-b-rel} b=-a^m/2^{m-1}.\eeq Now
consider the natural forgetful map
$$\Deff(f)\to \Deff(C)\times\Deff(X)$$
which, in terms of coordinates, is just given by projection to the $ab$ plane.
Its fibre is just $\Deff(f^!\t_X)$ which, as we have seen above, is smooth. Therefore this is
is a smooth morphism onto the curve $D$ with equation \refp{a-b-rel}., whose fibre coincides with $\Deff(f^!\t_X)$. This curve $D$ in turn projects smoothly to the $a$-axis, viz. $\Deff(C)$, but is ramified of degree $m$ over the $b$ axis, viz. $\Deff(X).$

\newsection{The projective case} The purpose of this section
is to adapt the construction of the tangent atom for an affinely
embedded scheme $X\to P$, as given in the previous section,
to the case of a projectively embedded scheme $X\to\P$. Whereas
in the affine case the basic idea was to take for tangent atom
$\t_X(P)$ the infinitesimal automorphisms of a resolution of
$X$, the idea here, not surprisingly, is to
construct the tangent atom $\t^.=\t_X(\P))$ by looking instead
 at graded resolutions of the homogeneous ideal and homogeneous
coordinate ring of $X$.\par
Unlike in the affine case, it is important here to work with the
sheafified version $\ttt^.$ of the tangent atom, in addition
to the graded module version $\t^.$. This has to do with the
fact that free graded modules are not acyclic. Indeed when it
comes to deformations, the non-acyclicity expresses itself in
the difference between the
appropriate cohomology of $\t^.$, which  yields the
\emph{projective} deformations, while the corresponding
(hyper-)cohomology of $\ttt^.$ yields all
(abstract) deformations.
In general, a scheme $X$ with $H^2(\O_X)\neq 0$ will admit
nonprojective deformations, and the associated classes in the
appropriate (Jacobi) cohomology of $\ttt_X$,
considered as formal power series, will be necessarily
\emph{non-convergent}.\par
 Although a tangent \sela $\t_{\bullet X}$, yielding all
deformations,  will be constructed in
the next section in the generality of arbitrary algebraic
schemes $X$, the sheafified tangent atom $\ttt^._X$
provides an adequate, and considerably simpler, substitute for
$\t_{\bullet X}$ in case $X$ is projective, as long as one
restricts to infinitesimal deformations, avoiding questions of
convergence.
\subsection{Differential operators on projective spaces}
Fix a projective space $\P=\P(V)$ and let $S=S_\P=\Sym^.(V^*)$
be the corresponding homogeneous coordinate ring. We will
generally identify a coherent sheaf $E$ on $\P$ with
its 'de-sheafification', i.e. graded $S$-module
$$\Gamma_*(E)=\bigoplus\limits_{n\in\Z} \Gamma(E(n)).$$
We will also identify a graded $S$-module $M$ with a graded free
resolution of $M$, as well as the sheafification
$\tilde{M}$.\par
Let
\beq\That_\P=\D^1(\O_\P(1), \O_\P(1))\eeq be the Lie algebra of
1st order
differential endomorphisms of $\O_\P(1)$ (which we may
identify with its module version, which is
the shift $S[1]$). As is well known, we have a
natural Lie isomorphism
\begin{equation}\label{dif}\That_\P\simeq \D^1(\O_\P(m),
\O_\P(m))
\end{equation}
 for any $m\neq 0$, and  a sheaf isomorphism
  $$\That_\P\simeq V\otimes\O_\P(1).$$ 
  The subsheaf
$$E=\D^0(\O_\P(1), \O_\P(1))\simeq \O_\P\subset \That_\P$$ is a
Lie ideal
 and the Euler
sequence \begin{equation}\label{euler}0\to E\stackrel{e}{\to}
\That_\P\to T_\P\to 0\end{equation} or in the module version
$$0\to S\to V\otimes S[1]\to \Gamma_*\That_\P\to 0$$
 is a sequence of Lie algebras and
homomorphisms, where $E$ is an abelian ideal and the action of
$T_\P$ on $E=\O_\P$ is the standard one. There is a natural
action pairing called the \emph{standard action}
$$\That_\P\times\O_\P(m)\to\O_\P(m).$$ For $m\neq 0$ this is
clear from
the identification (\ref{dif}), while for $m=0$ it comes via
(\ref{euler}) from the action of $T_\P$ on $\O_\P$. For
$m=0$, this action is a derivation in the usual sense. For
all m, it is a 'derivation relative to the action on $\O_\P$',
in the sense that
$$v(af)=av(f)+v(a)f, v\in\That_\P, a\in \O_\P,
f\in\O_\P(m).$$\par
From \refp{dif}. we deduce, via post-multiplication,
isomorphisms
\beq\label{diff-twist}\That_\P(k)\simeq \D^1(\O_\P(m),
\O_\P(m+k)), \forall m,k.\eeq This again defines the
\emph{standard action}
$$\That_\P(k)\times\O_\P(m)\to\O_\P(m+k).$$

We note that over each standard affine $P=D_{X_i}\subset\P$, the
identification $\O_P\stackrel{X_i}{\to}\O(1)|_P$ induces a
Lie-theoretic splitting of \refp{euler}., i.e.
\beq\label{eulersplit}\That_\P|_P\simeq E\oplus T_P.\eeq
 Moreover, the $E$ summand (i.e. the Euler operator) acts
trivially on $\O(m)|_P$
for all $m$.
\par \subsection{Construction}
Fix a closed subscheme $X\subset\P=\P(V),$ and consider a free
resolution of its ideal sheaf $\I=\I_X$:
\beq\label{res-proj}...\to \F^i=\bigoplus\F^i_\alpha\to...\to
\F^0\to\I\eeq
where we may assume  $$\F^i_\alpha=\O_\P(-m_{i\alpha}),
m_{i\alpha}>0, \forall i.$$ One way to obtain some
(though not all) sheaf resolutions \refp{res-proj}. is to start
with a graded resolution $F^.$ of the graded ideal $I^.=I^._X$
and sheafify. Resolutions obtained in this manner are said to be
\emph{full}. A general sheaf resolution \refp{res-proj}. is
obtained by sheafifying a complex of free modules which, at each
stage, is exact \emph{in sufficiently high degrees}, i.e. whose
cohomology is supported on the irrelevant prime of $S$.\par
As before, we augment $\F^.$ to a resolution $\F^._+$ of $\O_X$
by setting $\F^1_+=\O_\P$. For each $\alpha, \beta$, we denote
by \beq\label{delta}\delta^i_\ab\in\Hom(\F^i_\alpha,
\F^{i+1}_\beta)=\O_\P(m_{i\alpha}-m_{i+1, \beta})\eeq the
appropriate component of the differential of the resolution.

Then we may define the (sheafy) normal atom as
$$\tn=\tn_{X/\P}=Hom^.(\F^.,\F^._+).$$
This is a sheaf of Lie atoms (in fact, dgla's) on $\P$.
The restriction of $\tn_{X/\P}$ on a standard affine open
$P\subset\P$ is just the sheafified version of the affine normal
atom $\tn_{X\cap P/P}$. On the other hand, applying $\Gamma_*$
we get the global graded version over $S_\P$, denoted
$\n_{X/\P}$.\par Now the idea, as in the affine case, is to
construct the sheafy tangent atom $\ttt^.$ as the mapping cone
of a suitable map $$T_\P\otimes \F^._+\to \tn_{X/\P}$$ where
$T_\P$ itself is identified with a mapping cone as in the Euler
sequence \refp{euler}..
The global graded case is analogous.\par
Now the terms of $\ttt^.$ are defined as
follows
\begin{eqnarray}
\ttt^i&=&(\That_P\otimes \F^{i+1}_+)\oplus (E\otimes
\F^{i+2}_+)
\oplus \gl^i(\I),\ \  i\leq 0,\\
&=&\gl^i(\I)\oplus\Hom(\F^{-i+1},\O_\P),\ \ \ \ \   i\geq 1.
\end{eqnarray} where, as before, we identify
\beq\gl^i(\I)=\gl^i(\F^.)=\bigoplus\limits_j\Hom(\F^j,
\F^{i+j}).\eeq
So, schematically,\beq
\t_X:...\to\begin{matrix}\gl^{-1}(\I)&&\gl^0(\I)&&\gl^1(\I)\\
\oplus&\to&\oplus&\to&\oplus\\
(\That_P\otimes \F^0)\oplus E&&\That_P&&\check{\F}^0\end{matrix}
\to ...
\eeq (with the middle term in degree 0).
As in the affine case, $\ttt^.$ will be constructed so as to
admit $\tn$ as subcomplex and $T_\P\otimes \F^._+$ as quotient
complex, thus the only components of the differential that
require definition are those going from $T_\P\otimes \F^._+$ to
$\tn$.\par Those components are defined as
follows.\begin{itemize}
\item For $i=0$, the map $\That_\P\to\tn^1$ is
\beq v\mapsto
\bigoplus\limits_{j,\alpha,\beta}v(\delta^j_\ab)\eeq
(cf. \refp{dif}. and \refp{delta}.)
\item For $i<0$, the map $\That_\P\otimes
\F^{i+1}_\gamma\to\tn^{i+1}$ is
\beq v\otimes a\mapsto
\bigoplus\limits_{j,\alpha,\beta}\mu(a,v(\delta^j_\ab))\eeq
where $\mu$ (abusively) denotes the composite map
$$\F^k\otimes\Hom(\F^r, \F^s)\to Hom(\F^r,
\F^{k}\otimes\F^s)\stackrel{\mu}{\to}\Hom(\F^r, \F^{k+s-1})$$
where the latter $\mu$ is the map induced by the multiplication
on $\O_X$ as in \S3.
\item form $i=-1$, the map $E\to \tn^0$ is
\beq
1\mapsto\bigoplus\limits_{j,\alpha}m_{j\alpha}\id_{\F^j_\alpha}.\eeq
\item For $i<-1$, the map  $E\otimes\F^{i+2}\to \tn^{i+1}$ is
\beq 1\otimes a\mapsto
\bigoplus\limits_{j,\alpha}m_{i\alpha}\mu(a,\id_{\F^j_\alpha}).\eeq

\end{itemize}
All these formulas are dictated, on the one hand, by the natural
action of $\That_\P$ on the various $\F^i_\ab$, and on the other
hand, by the concept of defining the differential going into\nl
$\tn\subset\Hom(\F^._+, \F^._+)$ via commutator with the
differential of $\F^._+$, as in \refp{dTXP-def-bracket}. and
\refp{dTX-def-bracket}.. Also, via the splitting
\refp{eulersplit}., they are compatible with the definition of
the affine tangent atom $\t_{X\cap P}(P)$ for each standard
affine $P=D_{X_i}\subset\P$.
\par
Next we must show that $\ttt$
is a complex, i.e. that $d^2=0.$ On all terms except those
involving $E$, the proof is similar to the affine case, based on
the fact that $\That_\P$ acts as derivations. Alternatively, one
can argue based on \refp{dTXP-def-bracket}. and
\refp{dTX-def-bracket}..
  On the terms involving $E$, i.e. $E\otimes \F^i$, the required
vanishing is an immediate consequence of Euler's identity,
i.e. the fact that $E$ acts as $m.\id$ on $\O(m)$.\par
Thus we have defined $\ttt=\ttt_X(\P)$ as a complex of
 sheaves (of graded free $\O_\P$ modules) on $\P$. We note that
the restriction of $\ttt$ on a standard affine open
$P=D_{X_i}\subset \P$ is quasi-isomorphic to the sheafy version
of the tangent atom $\ttt_{X\cap P}(P)$.
 Note also that by construction, there is a subcomplex
$\ttt_{X/\P}$ of $\ttt$ (omitting all the $T_\P\otimes \F^.$
terms), such that $\ttt$ is the mapping cone of a map
$T_\P\otimes \F^.\to\ttt_{X/P}$.\par
 Next, we define an action of $\ttt_{X/\P} $ on the
 complex $\F^.(m)\sim \I(m), \forall m$, analogous to the affine
case.
 We let $\tn$ act in the obvious way, as in the affine case. We
let $\That_\P$ act in the standard way on each $\O(k)$ summand
of $\F^i(m)$ and similarly for $E$ (which acts on $\O(k)$ as
multiplication by $-k$).  We note that the restriction of this
action on a standard affine $P\subset\P$ is compatible with the
action of $\t_{X/P}$ on $F^.$ defined previously, via
 the splitting \refp{eulersplit}.
 with $E$ acting trivially on $\O_P$.\par
 Next, the action of $\ttt_{X/\P}$ on $\F^.$ extends to an
action of $\ttt_{X/\P}$ on $\F^._+(m)$, where $\tn$ acts
trivially on $\F^1_+(m)=\O_\P(m)$. Similarly, there is an
action of $\ttt_{X/P}$ on $\tn$, where $\tn$ acts on itself via
bracket (adjoint action), and $\That_\P$ and $E$ act on each
line bundle summand $\O(k)\subset \tn$ in the standard
way.\par
 Next, the action of $\ttt_{X/\P}$ on $\F^._+(m)$ extends to an
action of $\ttt$ on $\F^._+(m)$ via the standard recipe
\beq\label{act-TF} \act{a\otimes v}{b}=\mu(a, \act{v}{b}), a\in
\F^i, b\in\F^j_+, v\in\That_\P\eeq where $\act{v}{b}$ refers to
the pairing defined previously as part of the action of
$\t_{X/\P}$

 \par Next, we define a bracket on $\ttt$ and $\ttt_{X/\P}$:
this is defined as the usual bracket (signed commutator) on
 $\n$ and the standard bracket on $T_\P$ (which is compatible
with the bracket on $\That_\P$). As for mixed terms, we define
$$[v,a]=\act{v}{a},\ \  v\in \That_\P, a\in\tn$$ where again
$\act{v}{a}$ refers to the pairing defined previously.
This suffices to define a bracket on $\ttt_{X/\P}$.
 We extend
this to a bracket on $\ttt$ by setting
$$[a\otimes v, b]=\act{a\otimes v}{b}$$
cf. \refp{act-TF}. As before, the necessary axioms (e.g. Jacobi identity, action rule etc.) can either be easily verified directly or deduced from the corresponding properties on the
affine pieces $X\cap P\subset P$.\par
This completes the construction of the tangent dgla $\ttt_X(\P)$ for any projective scheme $X\subset \P$. As in the affine case, we can show that\begin{enumerate}
\item the dgla quasi-isomorphism class of $\ttt_X(\P)$ depends only on the subscheme $X\subset \P$ and not on the resolution and other choices;
\item there is a reduced version $\ttt_X^\red(\P)$, again independent of the resolution;
\item the weak equivalence classes of $\ttt_X(\P)$ and $\ttt_X^\red(\P)$ are independent of the embedding $X\to \P$ and depend on the isomorphism class of $X$ only.
\end{enumerate}
\begin{rem} If$X$ is as above (closed in $\P$) and $X'\subset X$ is any open subset, the restriction $\ttt_{X}(\P)|_{X'}$ yields
a dgla sheaf on $X'$ that will be called it tangent dgla and denoted $\ttt_{X'}(\P)$. This extends the notion of tangent dgal to the case of an arbitrary quasi-projective scheme.
\end{rem}
\newsection{Algebraic schemes: Tangent \sela}\label{algebraic}
The purpose of this section is to extend the notion of tangent
object (dgla
or atom) of affine schemes and their maps to the case of
general
(algebraic) schemes. This extension will take the form of a
semi-simplicial Lie algebra or \sela. The reason for this added
level of complexity is that, given an affine covering
$(X_\rho)$
of
a scheme, each tangent dgla $\t_{X_\rho}$ is in reality defined
with respect to some particular affine embedding $X_\rho\to
P_\rho$, and these embeddings are not mutually related.
Consequently, the $\t_{X_\rho}$ need not glue together to
a dgla, because the restrictions of $\t_{X_\rho}$ and
$\t_{X_\sigma}$ on $X_\rho\cap X_\sigma$ are only weakly
equivalent
(and indirectly so at that), and this relation is too weak for
ordinary gluing. Nonetheless the relation of weak equivalence
is
strong enough to yield a \sela.
\subsection{Construction}
 Here we construct the tangent \sela of a
separated algebraic scheme $X$ over $\C$ (the separatedness
does
not
seem to be essential, but is convenient). To fix ideas we focus
on the non-sheafy version,
though the (coherently) sheafified version can be constructed
in
the
same way. Let $(X_\rho)$ be an affine open covering of $X$
indexed
by a well-ordered set, and for each $\rho$ let $P_\rho$ be an
affine
space with a closed embedding
\begin{equation}\iota_\rho:X_\rho\subset P_\rho.\end{equation}
 Set $B_\rho=A_{P_\rho}$. We call the
system $(X_\rho\subset P_\rho)$ an \emph{affine embedding
system}
for $X$. Note that via $$X_\rho\cap X_\sigma=\Delta_X\cap
X_\rho\times X_\sigma\subset X\times X,$$ $X_\rho\cap X_\sigma$
embeds as a closed
subscheme of $X_\rho\times X_\sigma$, hence of $ P_\rho\times
P_\sigma.$ Similarly, for any multi-index $\rho_0<...<\rho_k$,
we
define
\begin{equation} X_{(\rho_0,...,\rho_k)}=
X_{\rho_0}\cap...\cap X_{\rho_k},\ \
P_{(\rho_0,...,\rho_k)}=P_{\rho_0}\times...\times P_{\rho_k}
\end{equation}
and the natural closed embedding
\begin{equation}\iota_{(\rho_0,...,\rho_k)}:\label{k-fold
inters}
X_{(\rho_0,...,\rho_k)}\subset
P_{(\rho_0,...,\rho_k)},\end{equation} and we denote the ideal
of
the latter by $I_{(\rho_0,...,\rho_k)}$. We call the system
\beq
{(} X_{(\rho_0,...,\rho_k)}\subset P_{(\rho_0,...,\rho_k)},
(\rho_0<...<\rho_k), k\geq 0{)}\eeq the \emph{simplicial
extension} of the affine embedding system $(X_\rho\subset
P_\rho)$.
Note that the defining equations for the image of
$\iota_{(\rho_0,...,\rho_k)}$ consist of defining equations for
the
images of individual embeddings $\iota_{\rho_i}$, together with
equations for the small diagonal on $X^{k+1}.$ The latter are
of
course generated by the pullbacks of the equations of the small
diagonal in $X^k$ via the various coordinate projections
$X^{k+1}\to
X^k$. Therefore, it is possible to choose mutually compatible
free
resolutions for all the $I_{(\rho_0,...,\rho_k)}$, and we
denote
these by $F^._{(\rho_0,...,\rho_k)}$. In fact, we may assume
that
\beq
F^1_{(\rho_0,...,\rho_k)}=B_{(\rho_0,...,\rho_k)}:=\bigotimes\limits_0^k
B_{\rho_i},\eeq \beq
F^i_{(\rho_0,...,\rho_k)}=\bigoplus\limits_{j=0}^k
(F^i_{\rho_j}\otimes
B_{(\rho_0,...,\rho_k)})\oplus\Delta^i_{(\rho_0,...,\rho_k)},
i\leq
0,\eeq where $\Delta^._{(\rho_0,...,\rho_k)}$ is a lifting to
$B_{(\rho_0,...,\rho_k)}$ of a free resolution of the small
diagonal\nl
$X_{(\rho_0,...,\rho_k)}\subset\prod\limits_{j=0}^kX_{\rho_j}$
and
moreover for any biplex
$$\rho^k=(\rho_0,..,\rho_k)\subset \rho^{k+1}=
(\rho_0,...,\rho_{k+1}),$$ if we let
$$\pi_{\rho^{k+1},\rho^{k}}:P_{\rho^{k+1}}\to P_{\rho^{k}}$$
denote
the natural projection, then we have a direct summand inclusion
\begin{equation}\label{Odiff}\pi_{\rho^{k+1},\rho^{k}}^*F^._{\rho^{k}}:=
\pi_{\rho^{k+1},\rho^{k}}\inv F^._{\rho^{k}}\otimes
B_{{\rho^{k+1}}}\to F^._{\rho^{k+1}}\end{equation} Putting
together
these groups and maps, and twisting by the appropriate sign,
i.e.
$\epsilon(\rho^k,\rho^{k+1})$, we get an 'extrinsic \v{C}ech
(double) complex'
\begin{equation} \xc(\O_X):
\bigoplus\limits_{\rho^0}F^._{+\rho^0}\to...\to
\bigoplus\limits_{\rho^k}F^._{+\rho^k}\to...\end{equation} A
similar
construction can be applied to any coherent sheaf on $X$.
Actually
this complex is quasi-isomorphic to the usual \v{C}ech complex
of
$ \O_X$, but we can do more with it. Note that the map
\refp{Odiff}.
gives rise to a direct summand inclusion
$$\pi_{\rho^{k+1},\rho^{k}}^*(\gl^.(I_{\rho^{k}}))\to
\gl^.(I_{\rho^{k+1}}),$$ whence a dgla map
$$\delta_{\rho^{k},\rho^{k+1}}:
\n_{X_{\rho^{k}}/ P_{\rho^{k}}}\to
\n_{X_{\rho^{k+1}}/P_{\rho^{k+1}}}$$ Then we can similarly
construct
a 'normal \sela ' $$\n_{X_\bullet/P_\bullet}:...\to
\bigoplus\limits_{\rho^k}\n_{X_{\rho^{k}}/
P_{\rho^{k}}}\to...$$
Likewise, we have an 'ambient tangent complex' $T_{P_\bullet}$
and
$T_{P_\bullet}\otimes\xc(\O_X)$ and a map
\begin{equation}\label{tangsela}
T_{P_\bullet}\otimes\xc(\O_X)\to
\n_{X_\bullet/P_\bullet}\end{equation} Finally, we define the
\emph{tangent \sela} of $X$ (with reference to the simplicial
system
$(X_\bullet, P_\bullet)$) to be the mapping cone of this, and
denote
it by $\t_{X\bullet}(X_\bullet, P_\bullet)$ or simply
$\t_{X\bullet}$. This is the \sela whose value on the simplex
$\rho^k$ is the dgla $\t_{X_{\rho^k}}(P_{\rho^k})$. By Theorem
\ref{tang sela}, there is an associated Jacobi-Bernoulli
complex
$J^.(\t_{X_\bullet})$, which we denote by $J^._X$ and refer to
as
the \emph{Jacobi-Bernoulli complex of} $X$. Up to filtered,
comultiplicative weak equivalence, it depends only on the
isomorphism class of $X$ as scheme over $\C$. Therefore the
\emph{deformation ring} of $X$ $$R_X=\C\oplus\HH^0(J^._X)^*$$
is
canonically defined. In the next section we relate $R_X$ to
flat
deformations of $X$ over artin rings. For any artin local
$\C$-algebra $S$, we set
$$J^._{X,S}=J^.(\t_{X_\bullet}\otimes\m_S)$$ and note that via
the
natural map $J^._{X,S}\to J^._X\otimes\m_S$, any class
$\epsilon\in\HH^0(J^._{X,S})$ yields a local homomorphism
('classifying map')
$$^t\epsilon: R_X\to S.$$
\par As in the affine case, this construction may be extended
from the case of schemes to that of maps. Thus let
$$f:X\to Y$$ be a morphism of schemes. Then we can choose
respective
affine coverings
$$X_\alpha\to P_\alpha, Y_\alpha\to Q_\alpha \textrm{\ \ such
that\
\ } f(X_\alpha)\subset Y_\alpha.$$ Then for each simplex
$\rho^k$,
the restriction of $f$ yields a morphism
$$f_{\rho^k}:X_{\rho^k}\to Y_{\rho^k},$$
and for this we have an associated tangent dgla
$\t_{f_{\rho^k}}$.
Putting these together, we get a tangent \sela (with respect to
the
given affine coverings)
$$\t_f:...\to \t_{f_{\rho^k}}\to...$$ As before, $\t_f$ can be
described as
a mapping cone:
\begin{equation}\label{tf-cone}\t_f={\mathrm{cone}}(
\t_{X\bullet}\oplus\t_{Y\bullet}\to
f^!\t_{Y\bullet}).\end{equation}
Thus we have \sela morphisms $$\t_f\to \t_Y, \t_f\to \t_X,
f^!\t_Y[1]\to\t_f.$$ Correspondingly, we have a
Jacobi-Bernoulli
complex $J^._f$, a deformation ring $R_f$ together with maps
$R_X\to
R_f, R_Y\to R_f.$\begin{rem}{\rm When $X$ is smooth, its
tangent
\sela is equivalent to a dgla, e.g. the Kodaira-Spencer
algebra,
a
soft dgla resolution of the tangent sheaf.}\end{rem}
\begin{rem} When $X$ is quasi-projective, i.e. a locally
closed subscheme of a projective space $\P$ with closure $\bar{X}$, we have constructed in the previous section
a tangent dgla sheaf $\ttt_{\bar{X}}(\P)$
\end{rem}
\subsection{Basic properties} We summarize the basic properties
of
the tangent \sela $\t_{\X_\bullet}(\X_\bullet,P_\bullet)$ as
follows.\begin{thm} \begin{enumerate}\item
$\t_{\X_\bullet}(\X_\bullet,P_\bullet)$ is a \sela and
$\xc(\O_X)$-module and acts on $\xc(\O_X)$.\item The weak
equivalence class $\t_{X_\bullet}$ of
$\t_{\X_\bullet}(\X_\bullet,P_\bullet)$ depends only on $X$ and
is
functorial.\item If $X$ is smooth, $\t_{X\bullet}$ is weakly
equivalent to the (\v{C}ech complex of) the usual tangent
algebra
$T_X$.
\end{enumerate}
\end{thm}\begin{proof}(i) The fact that
$\t_{\X_\bullet}(\X_\bullet,P_\bullet)$ is a \sela is clear
from
the
construction, and the actions by and on $\xc(\O_X)$ are also
clear
component by component.\par (ii) Any two affine embedding
systems
$(X_\bullet, P_\bullet), (X_\bullet, Q_\bullet)$ are dominated
by a
third, say $(X_\bullet, R_\bullet)$, and we get maps
$$\t_{X_\bullet}(X_\bullet, P_\bullet)\to
\t_{X_\bullet}(X_\bullet,
R_\bullet)\leftarrow\t_{X_\bullet}(X_\bullet, Q_\bullet).$$
These
maps are termwise weak equivalence hance, by a standard
spectral
sequence argument, overall weak equivalences as well. This
proves
the independence of $\t_{X_\bullet}(X_\bullet, P_\bullet)$ on
the
affine embedding system used to define it. The functoriality is
proved similarly.\par (iii) This follows from Corollary
\ref{smoothaffine} above, since in he smooth case we may take
for
each $\t_{X_\rho}$ the usual tangent algebra, and these glue
together properly on overlaps.
 \end{proof}
\par\newsection{Deformations of schemes}
\subsection{Classification}
Deformations of an algebraic
scheme $X/\C$ can be classified in terms of the associated
tangent
\sela $\t_X$ and its Jacobi-Bernoulli cohomology. Consider
first
the
case of an affine scheme $X\subset P$ (notations as in
\S\ref{affine}). Let
$S$ be a local artinian $\C$-algebra. Then a flat deformation
of
$X$
over $S$ is determined by, and determines, up to certain
choices, an
element
\begin{equation}\label{aff_ks_elt}\phi\in\t^1_X(P)\otimes\m_S=
\Hom^1(F^._X,
F^._{+X})\otimes\m_S\end{equation} known as a
\emph{Kodaira-Spencer
cochain}, which satisfies the integrability condition
\begin{equation}\label{aff_integ}\partial\phi=-\half[\phi,\phi].
\end{equation}
The deformation corresponding to $\phi$ can be determined e.g.
as
the subscheme of $P\times\Spec(S)$ having $(F^._X\otimes S,
\partial +\phi)$ as resolution; we may denote this by
$X^\phi$.\par
Now globally, let $X$ be an algebraic scheme over $\C$ and as
in
\S\ref{algebraic} choose an affine embedding system
$$\iota_\rho:X_\rho\to
P_\rho.$$ This gives rise as in \S\ref{algebraic} to a
representative for the
tangent \sela $\t_X$. Now suppose given a deformation $\mathfrak
X$
of
$X$ over $S$ as above. This restricts for each $\rho$ to a
deformation $\mathfrak X_\rho$ of $X_\rho$, whence a
Kodaira-Spencer
cochain $$\phi_\rho\in\t^1_{X_\rho}(P_\rho)\subset
\t^1_X(P_\bullet),$$ satisfying an integrability condition as
in
(\ref{aff_integ}), so that the restricted deformation
$\mathfrak
X_\rho$ is $S$-isomorphic to  $X_\rho^{\phi_\rho}$. Moreover,
the
fact that $\phi_\rho$ and $\phi_\sigma$ restrict to equivalent
deformations of $X_{\rho\sigma}\subset P_{\rho\sigma}$ yields
an
$S$- isomorphism
\begin{equation}\label{sigma_iso_rho}X_\rho^{\phi_\rho}\cap
X_\sigma\simeq X_\sigma^{\phi_\sigma}\cap X_\rho;\end{equation}
both
of these are closed subschemes of
$P_{\rho\sigma}\times\Spec(S)$
and
the isomorphism (\ref{sigma_iso_rho}) extends to an $S$-
automorphism of $P_{\rho\sigma}\times\Spec(S)$, necessarily of
the
form
$$\exp(t_{\rho\sigma}), t_{\rho\sigma}\in
T_{P_{\rho\sigma}}\otimes\m_S.$$ Then we get two resolutions of
$X_\rho^{\phi_\rho}\cap X_\sigma$, the 'original' one with
differential $\partial +\phi_\rho$, and the one pulled back
from
$X_\sigma^{\phi_\sigma}\cap X_\rho$, whose differential is
$\partial
+\phi_\sigma.$ It is easy to see and well known that the two
resolutions differ by an isomorphism of the form $\exp(u_\rs)$
where
$u_\rs\in\gl^0(F^._{X_\rho})\otimes\m_S.$ Thus all in all there
is a
(uniquely determined) element
$$\psi_{\rho\sigma}\in\t^0_{X_{\rho\sigma}}\otimes\m_S=
(T_{P_\rs}\oplus\gl^0(F^._{X_\rs}))\otimes\m_S$$ such that
\begin{equation*}
\exp(\psi_{\rho\sigma})(\partial+\phi_\sigma)
\exp(-\psi_{\rho\sigma})=\partial+\phi_\rho.\end{equation*}
More
explicitly, the latter equation reads $$\partial
\exp(\psi_{\rho\sigma})=-\phi_\rho\exp(\psi_{\rho\sigma})
+\exp(\psi_{\rho\sigma})\phi_\sigma.$$ Multiplying by
$\exp(-\psi_{\rho\sigma})$, we get
$$D(-\ad(\psi_{\rho\sigma}))(\partial\psi)=
-\exp(-\psi_{\rho\sigma})\phi_\rho\exp(\psi_{\rho\sigma})+\phi_\sigma$$
where $D(x)$ is as in (\ref{cd}), with inverse $C(x)$ the
generating
function for the Bernoulli numbers. Multiplying the latter
equation
by $C(-\ad(\psi_{\rho\sigma}))$ and using the simple identity
$C(-x)e^{-x}=C(x)$, we conclude
\begin{equation}\label{phipsi_cocycle}\partial\psi_{\rho\sigma}=
-C(\ad(\psi_{\rho\sigma}))\phi_\rho+C(-\ad(\psi_{\rho\sigma}))\phi_\sigma.
\end{equation}
The latter equation (\ref{phipsi_cocycle}) is just the
condition
that the component in bidegree $(0,1)$ of   the coboundary
$\partial
\epsilon(\phi_\bullet, \psi_\bullet)$ should vanish. Also, by
construction, we clearly have
\begin{equation}\label{psi_cocycle}
\exp(\psi_{\rho\sigma})\exp(\psi_{\sigma\tau})=
\exp(\psi_{\rho\tau})\end{equation} Thus, putting
(\ref{aff_integ},
\ref{phipsi_cocycle}, \ref{psi_cocycle}) together and using the
derivation property of the differentials in the
Jacobi-Bernoulli
complex,
 $\epsilon(\phi_\bullet, \psi_\bullet)$ is a special
multiplicative
cocycle in the Jacobi-Bernoulli complex $J(\t_X)\otimes \m_S.$
Conversely, given a special multiplicative cocycle
$\epsilon(\phi_\bullet, \psi_\bullet)$ with values in $S$, the
$\phi_\bullet$ data yields a collection of deformations of the
affine pieces of $X$, while the $\psi_\bullet$ glues these
deformations together. These processes are inverse to each
other
up
to an automorphism, and are precise mutual inverses when there
are
no automorphisms. Hence \begin{thm} Let $X$ be an algebraic
scheme
over $\C$ such that $H^0(\t_X)=0.$ Then for any local artin
$\C$-algebra $S$, there is a natural bijection between the set of
equivalence classes of flat deformations of $X$ over $S$ and
the
set
of local homomorphisms from $R_X$ to $S$.\end{thm}
\subsection{Obstructions} Let $S$ be a local artin algebra,
$I<S$
an ideal contained in the socle $\mathrm{ann}_S(\m_S)$ and
$\bar{S}=S/I.$ Let $\bar{\epsilon}=\epsilon(\bar{\phi}_\bullet,
\bar{\psi}_\bullet)$ be a special multiplicative cocycle with
coefficients in $\bar{S}$. Let $\phi_\bullet, \psi_\bullet$ be
arbitrary liftings of $\phi_\bullet, \psi_\bullet$ with
coefficients
in $S$. Thus, $\epsilon=\epsilon(\phi_\bullet, \psi_\bullet)$
is
not
necessarily a cocycle. However, it is easy to check that the
coboundary $\partial \epsilon$ lies in  $$(F_1J_X)^1\otimes
I=(K^0(\t_X)^2\oplus K^1(\t_X)^1\oplus K^2(\t_X)^0)\otimes I$$
(this is because it dies modulo $I$) and moreover, that
$\partial\epsilon$
is
a cocycle for $\mathrm{tot}(K^.(\t_X)^.)\otimes I$ (the latter because it
is a
cocycle for $J^._{X,S}$). Thus, we obtain a cohomology class
\begin{equation}\label{obstruction}
\mathrm{ob}(\bar{\phi}_\bullet,
\bar{\psi}_\bullet)\in H^2(\t_X)\otimes
I=H^2(\mathrm{tot}(K^.(\t_{X})^.)\otimes I.\end{equation} This
class
is independent of choices and represents the obstruction to
lifting
$\bar{\epsilon}$ to a special multiplicative cocycle with
coefficients in $S$.
\begin{example} Let $X$ be a quadric cone in $P=\P^3,$
 with equation $q=x_1x_2-x_3^2.$ Let\nl $f:L\to X$ be the inclusion of the line with equations $x_1=x_2=x_3$. Then $\t_X$ is as in Example \ref{quadcone} and $f^!\t_X[-1]$ is the Lie atom (special case of \sela)
$$I_L\t_X\to \t_X$$
This atom, being the mapping cone of $\t_f\to\t_X\oplus\t_Y$, classifies deformations of the map $f$ together with trivializations of the corresponding deformations of $X$ and $Y$.
Also, as $L$ is smooth we may identify $\t_L$ with the utual tangent algebra $T_L$, and there is a natural map
$$\t_L\to f^!\t_X$$
whose associated Lie atom, i.e. mapping cone, is $\n_{L/X}$, which classifies deformations of $L$ in $X$.
Then it is well known that
$$T_P\otimes \O_L\simeq T_L\oplus 2\O_L(1), N_{L/P}\simeq 2\O_L(1).$$
The map $T_P\otimes \O_L\to \O_L(2)$ vanishes on $T_L$ and has cokernel a skyscaper $\C$ at the vertex of the cone. The kernel of the induces map $N_{L/P}\to\O_L(2)$ is generated ny $x_0v, x_1v$ where
$$v=\partial/\partial x_1-\partial/\partial x_2\in T_P(-1)\subset \t_X^0(-1).$$
Thus $$\HH^0(\n_{L/X})\simeq\{a_0x_0v+a_1x_1v:a_0. a_1\in\C\}\simeq\C^2, \HH^1(\n_{L/X})=[x_0^2]\simeq\C.$$
Set $$g=x_0(x_2-x_1)\in I_L(2)=I_L\t_X^1.$$
Then the fact that $g=x_0v(f)$ shows that the element $\epsilon=\epsilon(x_0v, g)$ is a special multiplicative cocycle for $f^!\t_X$ with coefficients in $\C=\C[t]/(t)$. The obstruction to lifting $\epsilon$ over $S=\C[t]/(t^2)$ is precisely
$$x_0v(g)=-2x_0^2\neq 0.$$  Therefore $\epsilon$ does not lift.
In fact, since the first obstruction map goes onto $\HH^1(\n_{L/X})$, there are no higher obstructions.  Therefore in dual coordinates, the formal deformation ring is
$$\hat{R}=\hat{R}(\n_{L/X})\simeq\C[[a_0^*, a_1^*]]/((a_0^*)^2),$$
i.e. a double 'line';
the universal deformation over $\hat{R}$ is given by $\epsilon(a_0^*g, a_0^*x_0v+a_1^*x_1v)$. Since set-theoretically the locus of lines in $X$ is a conic, this shows the conic occurs with multiplicity 2.\par
A more involved calculation that we won't do in detail actually yields the equations of this double conic in a standard $\A^4$ neighborhood on $L$ in the Grassmannian of lines in $\P^3$, which is coordinatized via the vector fields $$x_0v, x_1v, x_0w, x_1w, w:=\partial/\partial x_1-\partial/\partial x_3.$$
 \end{example}
\subsection{Applications: relative obstructions, stable
subschemes and surjections} Let $f:X\to Y$ be an embedding of a
closed subscheme. We consider the question of 'relative
obstructions' i.e. obstructions to lifting a given deformation
of
$Y$ to a deformation of $f$. These are obstructions for the
mapping
cone of $\t_f\to \t_Y$, which is the same as that of $\t_X\to
f^!\t_Y$ (cf. \refp{tf-cone}.). Locally, if
 $Y\to P$ is an affine embedding, such obstructions have values
in
$\mathrm{Ext}^1(I_{X/Y}, \O_X)$, where
$I_{X/Y}=I_{X/P}/I_{Y/P}$
is
\emph{independent of $P$} even as complex up to
quasi-isomorphism
(not just weak equivalence). Hence, global obstructions also
have
values in $\mathrm{Ext}^1(I_{X/Y}, \O_X)$. In reasonably good
cases
though, the obstruction group can be narrowed considerably. The
following result sharpens Thm. 1.1 of \cite{Rhol}\begin{thm}
\label{relobst}
Let
$X\subset Y$ be a closed subscheme having no component
contained
in
the singular locus of $Y$. Then obstructions to deforming $X\to
Y$
relative to deforming $Y$ are in
\begin{equation}\label{image}\mathrm{Im(Ext}^1_Y(I_{X/Y}/I_{X/Y}^2,
\O_X)\to\mathrm{Ext}^1_Y(I_{X/Y},
\O_X)).\end{equation}\end{thm}
\begin{cor}Assumptions as in Theorem \ref{relobst}, assume additionally
that $X\to Y$ is a
regular embedding with normal bundle $N$ and $H^1(X,N)=0$. Then
$X\to Y$ is 'relatively unobstructed' or 'stable' relative to
$Y$, in the sense that it
deforms with every deformation of $Y$; furthermore the
Hilbert scheme of $Y$ is smooth of dimension $h^0(X,N)$
 at the point corresponding to
$X$.\end{cor}
This Corollary generalizes a result of Kodaira \cite{K1} in the smooth
case.\par To sketch the proof of Thorem \ref{relobst},
first working locally, let $F^._Y$ be
a
free resolution of $I_{Y/P}$ and extend it to a free resolution
$F^._X$ of $I_{X/P}$, such that, termwise,
$$F^i_X=F^i_Y\oplus F^i_{X/Y}$$ where $F^i_{X/Y}$ is a suitable
free
complement, $F^._Y\to F^._X$ is a map of complexes (though not
a
direct summand inclusion), and $F^._{X/Y}\otimes\O_Y$, which is
a
quotient complex of $F^._X$, is a free resolution of $I_{X/Y}$.
We
may also assume that $F^._{X/Y}$ contains a subcomplex
$F^._{X^2}$
(with termwise direct summands) resolving $I_X^2$. A
deformation
of
$Y$ yields a linear map
$$v:F^0_Y\to \O_X.$$ The obstruction to lifting this to a
deformation of $X$ is given by $$v\circ\delta:F^{-1}_{X/Y}\to
\O_X$$
where $\delta:F^{-1}_{X/Y}\to F^0_Y$ is the 'connecting map'
from
the resolution. $\delta$ takes a relation among generators of
$I_X\mod I_Y$ to the appropriate linear combination of
generators of
$I_Y$. Now, and this is the point, our assumption about
singularities means that no generator of $I_Y$ can be in
$I_X^2$, or
more precisely, that
$$I_{X}^2\cap I_Y\subset I_XI_Y.$$ This implies that
$\delta(F^{-1}_{X^2})\subset I_XF^0_Y$. Since $v$ is
$\O_P$-linear,
it follows that the composition $v\circ\delta$ is zero on
$F^{-1}_{X^2}$ and lives in the image as in (\ref{image}). This
proves the result in the affine case, and the extension to the
global case is straightforward.\qed.\par
If $Y$ is a projective locally complete intersection,
the reduced tangent algebra $\tred_Y$
has
an important \emph{numerical} property related to the dualizing sheaf of $Y$.
  Namely, viewing $\tred_Y$ as an element in the Grothendieck group $K(Y)$ of vector bundles (or coherent sheaves), it has a well-defined determinant line bundle, and we have
\beq\det(\tred_Y)=\omega_Y\inv=\omega_P\inv\det(N_{Y/P})\eeq
where $P$ is an ambient projective space.
Indeed this follows immediately from the construction of the dualizing sheaf and the discussion of Serre duality theory in
\cite{hart}. This, together with the fundamental dimension
inequality, Corollary \ref{dim-ineq}, can be applied to study deformations
of maps of curves to varieties:
\begin{thm}Let $f:X\to Y$ be a map of a smooth projective curve of genus $g$ to a projective locally complete intersection  variety of dimension $n$
such that $f(X)$ is not contained in the singular locus
of $Y$. Then the space $\Deff(f)$ of deformations of $f$ fixing $X,Y$ is of dimension at least
$\chi=X.\omega_Y\inv+n(1-g)$; moreover, L if $\chi>3, g=0, >1, g=1, or >0, g>1$, $f$ admits deformations so that $f(X)$ moves.
\end{thm}
\begin{proof} In our case $f^!\tred_Y$ is a finite
complex of vector bundles on $X$, acyclicl in degrees
$\neq 0,1$. By Riemann-Roch, we have
\beql{}{\chi(f^!\tred_Y)=(h^0-h^1)(f^!\tred_Y)=X.\omega_Y\inv+n(1-g).}
Therefore by Corollary \ref{dim-ineq}, our conclusion follows.
\end{proof}
\par
Next we consider an
application to surjections is (compare \cite{Rhol}, Thm
2.1):\begin{thm} Let $f:X\to Y$ be a projective morphism with
$$f_*(\O_X)=\O_Y, R^1f_*(\O_X)=0.$$ Then $f$ deforms with every
deformation of $X$.\end{thm}\begin{proof} Here we use the fact
that
the mapping cone of $\t_f\to\t_X$ is equivalent to that of
$\t_Y\to
f^!\t_Y.$ To prove the result it suffices to show the
following\newline (*)\emph{ the natural map $\t_Y\to f^!\t_Y$
induces a surjection on $H^1$ and an injection on $H^2.$}\par
Indeed
the $H^1$-surjectivity property implies that any first-order
deformation of $X$ lifts to a deformation of $f$; then the
$H^2$-injectivity property allows us to extend this
inductively,
via
obstruction theory (\S3.2), to $n$th order deformations. To
prove
(*), note that $f^!\t_Y$ can be represented by the tensor
product
$\t_Y\otimes\xc(\O_X)$, which is a double complex with terms
$\t^i_Y\otimes \xc^j(\O_X)$. Then the spectral sequence of a
double
complex (or an elementary substitute) yields our
conclusion.\end{proof} Morphisms $f$ satisfying the hypotheses
of
the Theorem occur in diverse situations, e.g. regular fibre
spaces
and resolutions of rational singularities. The Theorem says
that
those schemes which admit a structure such as $f$ form an open
neighborhood of $X$ in Moduli.
\begin{thm} Let $f:X\to Y$ be a proper surjective morphism
\'etale
in codimension 1 where $X$ is normal. Then any deformation of
$f$ is
determined by the associated deformation of $Y$.
\end{thm}\begin{proof} It will suffice to prove that the
mapping
cone $\t_X\to f^!\t_Y$ is exact in degrees $\leq 1$. Working
locally, we may view $X$ as a subscheme of $Y\times R$, $R$ an
affine space, and consider a free resolution $J^.$ of the ideal
of
$X$ in $Y\times R$. Let $K$ be the kernel of the natural map
$\hom(J^0,\O_X)\to\hom(J^{-1}, \O_X).$ Then $K$ is a
torsion-free
$\O_X$ module and the natural map $k:T_R\otimes\O_X\to K$ is an
isomorphism in codimension 1 by our assumption that $f$ is
\'etale
in codimension 1. As $X$ is normal, this easily implies $k$ is
an
isomorphism, which proves our assertion.
\end{proof}\begin{cor} A small resolution of a singularity $Y$
is
(infinitesimally) rigid relative to $Y$; hence the set of small
resolutions of a given $Y$ is at most countable.\end{cor}
Smallness is of course an essential hypothesis here. One
wonders,
however, whether 'countable' may be replaced by 'finite' in the
conclusion.\par Under stronger hypotheses on the size of the
exceptional locus, we can actually identify the deformations of
$f$
and $Y$ (compare \cite{Rhol}, Thm. 3.5):\begin{thm} Let $f:X\to
Y$
be a proper surjective morphism \'etale in codimension 2 with
$X$ locally
Cohen-Macaulay everywhere and lci in codimension 2.
Then any deformation of $Y$ lifts to a deformation of
$f$.\end{thm}\begin{proof}
We work as in the proof of Theorem \ref{compare-to-ext}. As there, consider the complexes
$$K_Y:...\F\inv_{1, Y}\otimes\O_Y\to F^0_Y\otimes\O_Y\to\Omega_Q\otimes\O_Y,$$
$$K_X:...\F\inv_{1, Y}\otimes\O_X\to F^0_X\otimes\O_X\to\Omega_P\otimes\O_X,$$
as well as the complex on $X$:
$$f^!K_Y:...\F\inv_{1, Y}\otimes\O_X\to F^0_Y\otimes\O_X\to\Omega_Q\otimes\O_X$$
As in the aforementioned proof, we have
$$\Hom^.(K_Y, \O_Y)\simeq\t_Y, \Hom^.(K_X, \O_X)\simeq \t_X,
\Hom^.(f^!K_Y, \O_X)\simeq f^!\t_Y.$$
Now let $K_f$ be the mapping cone of the functorial map $f^!K_Y\to K_X$ (cf. subsubsection \ref{functoriality};
 this cone is a complex in degrees $\leq 1$,
and its cohomology is supported on the non-\'etale
 locus of $f$, which has codimension $\geq 3$ in $X$. Thefore it follows from Ischebeck's theorem as in the proof of Theorem
 \ref{compare-to-ext} that we have
 $$\ext^i(K_f, \O_X)=0, i\leq 2.$$
 From this is follows that the natural map
 $$H^i(\t_X)\to H^i(f^!\t_Y)$$
 is bijective for $i\leq 1$ and injective for $i=2$;
  in other words, \par\noindent
 ${(\star)}$\ \ \ \ \ {\emph{the\ map\ }} $\t_X\to f^!\t_Y$
{\emph{\  is\ a\ weak\ equivalence.}}\par
   From this last statement the Theorem follows
   by a purely formal, standard argument that runs as follows.
Given an artinian local algebra $S$ and a local homomorphism $$\phi_Y:R(\t_Y)\to S,$$
we see by using $(\star)$ and artinian induction that the induced map
 $f^!\phi_Y: R(f^!\t_Y)\to S$ lifts to a map
 $\phi_X:R(\t_X)\to S$,  and this means precisely that $\phi_Y$ lifts to a deformation $(\phi_X, \phi_Y)$ of $f$.
\end{proof}
\begin{rem}The Theorem does not extend to morphisms with a
codimension-2 exceptional locus; for example, it fails for
a small resolution of a
3-fold ODP.
\end{rem}
\bibliographystyle{amsplain}
\bibliography{mybib}

\providecommand{\bysame}{\leavevmode\hbox to3em{\hrulefill}\thinspace}
\providecommand{\MR}{\relax\ifhmode\unskip\space\fi MR }
\providecommand{\MRhref}[2]{%
  \href{http://www.ams.org/mathscinet-getitem?mr=#1}{#2}
}
\providecommand{\href}[2]{#2}
\begin{thebibliography}{10}

\bibitem{hart}
R.~Hartshorne, \emph{Algebraic {G}eometry}, Springer, Berlin, 1977.

\bibitem{isch}
F.~Ischebeck, \emph{Eine dualit{\"a}t zwischen den funktoren {E}xt und {T}or},
  J. Alg. \textbf{11} (1969), 510--531.

\bibitem{HM}
J.Harris and D.~Mumford, \emph{The kodaira dimension of the moduli space of
  curves}, Inventiones math. (1984).

\bibitem{K1}
K.~Kodaira, \emph{A theorem of completeness of characteristic systems for
  analytic families of compact submanifolds of complex mainifolds.}, Ann. Math.
  \textbf{75} (1962), 146--162.

\bibitem{K2}
\bysame, \emph{Complex manifolds and deformations of complex structures},
  Springer, Berlin, 1986.

\bibitem{lich}
S.~Lichtenbaum and M.~Sclessinger, \emph{The cotangent complex of a morphism.},
  Trans. AMS \textbf{128} (1967), 41--70.

\bibitem{mat}
H.~Matsumura, \emph{Commutative algebra 2e}, Benjamin/Cummings, {Reading,
  Mass}, 1980.

\bibitem{sing}
Z.~Ran, \emph{Enumerative geometry of families of singular plane curves},
  Invent. math \textbf{97} (1989), 447.

\bibitem{Rhol}
\bysame, \emph{Stability of certain holomorphic maps}, J. Diff. Geom.
  \textbf{34} (1991), 37--47.

\bibitem{cid}
\bysame, \emph{Canonical infinitesimal deformations}, J. Algebraic Geometry
  \textbf{9} (2000), 43.

\bibitem{atom}
\bysame, \emph{Lie atoms and their deformation theory}, Geometric and
  Functional Analysis \textbf{18} (2008), 184--221.

\bibitem{Sernesi}
E.~Sernesi, \emph{Deformations of algebraic schemes}, Springer International,
  2006.

\bibitem{vara}
V.S. Varadarajan, \emph{Lie groups, {L}ie algebras and their representations},
  Springer, Berlin, 1973.

\end{thebibliography}
\end{document}